%% file: manuscript.tex
\documentclass[final,3p,times,fleqn]{elsarticle}

\usepackage{amsmath,amssymb,amscd,mathtools}
\usepackage{enumitem}
\usepackage{graphicx}
\usepackage{xcolor}
\usepackage{algorithm}
\usepackage{algpseudocode}
\usepackage{diagbox}
\usepackage{url}
\usepackage{tikz}
\usepackage{xparse}
\usepackage{booktabs}

\usepackage[colorlinks=true, allcolors=blue]{hyperref}
\usepackage{hypernat}

\usepackage{subcaption}
\usepackage{pgfplots}
\usepackage{pgfplotstable}
\biboptions{sort&compress}
\usetikzlibrary{matrix,backgrounds,shapes,decorations.pathreplacing,calligraphy,math}
\pgfdeclarelayer{myback}
\pgfsetlayers{myback,background,main}

\tikzset{mycolor/.style = {line width=1bp,color=#1}}%
\tikzset{myfillcolor/.style = {draw,fill=#1}}%

\NewDocumentCommand{\highlight}{O{blue!40} m m}{%
\draw[mycolor=#1] (#2.north west)rectangle (#3.south east);
}

\NewDocumentCommand{\fhighlight}{O{blue!40} m m}{%
\draw[myfillcolor=#1] (#2.north west)rectangle (#3.south east);
}

\newcommand{\wetting}{{w}}
\newcommand{\nonWetting}{{nw}}
\newcommand{\well}{{w}}
\newcommand{\res}{{r}}
\newcommand{\grid}{{\mathcal{T}}}
\newcommand{\faces}{{\mathcal{E}}}

\def\blkInterpolate{{\blkMat{P}_{\ell+1}^\ell}}
\def\blkRestrict{{\blkMat{R}_\ell^{\ell+1}}}

\def\blkInject{{\blkMat{Q}_\ell^{\ell+1}}}

\newcommand{\interpolate}[2][]{ \big(\Mat{P}_{{#2}^{#1}}\big)_{\ell+1}^\ell }
\newcommand{\restrict}[2][]{ \big(\Mat{R}_{{#2}^{#1}}\big)_\ell^{\ell+1}}

\definecolor{butter1}{rgb}{0.988,0.914,0.310}
\definecolor{chocolate1}{rgb}{0.914,0.725,0.431}
\definecolor{chameleon1}{rgb}{0.541,0.886,0.204}
\definecolor{skyblue1}{rgb}{0.247,0.524,0.912}
\definecolor{applegreen}{rgb}{0.55, 0.71, 0.0}
\definecolor{blue-green}{rgb}{0.0, 0.87, 0.87}
\definecolor{plum1}{rgb}{0.678,0.498,0.659}
\definecolor{scarletred1}{rgb}{0.937,0.161,0.161}

\newcommand{\tensorOne}[1]{\boldsymbol{#1}}

\def\TimeDomain{{\mathbb{T}}}

\renewcommand{\Vec}[1]{%
  \ifcat\noexpand#1\relax 
    \boldsymbol{#1}
  \else
    \mathbf{#1}
  \fi
}
\newcommand{\Mat}[1]{#1}

\newcommand{\blkVec}[1]{\Vec{#1}}
\newcommand{\blkMat}[1]{\Vec{#1}}

\newtheorem{rmk}{Remark}

\journal{Computers \& Mathematics with Applications}

\begin{document}
\title{Multilevel well modeling in aggregation-based nonlinear multigrid for multiphase flow in porous media}
\begin{frontmatter}
\author[casc]{Chak Shing Lee\corref{cor1}}
\cortext[cor1]{Corresponding author.}
\ead{cslee@llnl.gov}
\author[total]{Fran\c cois P. Hamon}
\ead{francois.hamon@totalenergies.com}
\author[aeed]{Nicola Castelletto}
\ead{castelletto1@llnl.gov}
\author[casc,psu]{Panayot S. Vassilevski}
\ead{vassilevski1@llnl.gov, panayot@pdx.edu}
\author[aeed]{Joshua A. White}
\ead{white230@llnl.gov}
\address[casc]{Center for Applied Scientific Computing, Lawrence Livermore National Laboratory, Livermore, CA 94550, USA}
\address[total]{TotalEnergies E\&P Research and Technology, Houston, TX 77002, USA}
\address[aeed]{Atmospheric, Earth, and Energy Division, Lawrence Livermore National Laboratory, Livermore, CA 94550, USA}
\address[psu]{Fariborz Maseeh Department of Mathematics and Statistics, Portland State University, Portland, OR 97201, USA} 




\begin{abstract}
A full approximation scheme (FAS) nonlinear multigrid solver for two-phase flow and transport problems driven by wells with multiple perforations is developed.
It is an extension to our previous work on FAS solvers for diffusion and transport problems.
The solver is applicable to discrete problems defined on unstructured grids as the coarsening algorithm is aggregation-based and algebraic.
To construct coarse basis that can better capture the radial flow near wells, coarse grids in which perforated well cells are not near the coarse-element interface are desired.
This is achieved by an aggregation algorithm proposed in this paper that makes use of the location of well cells in the cell-connectivity graph.
Numerical examples in which the FAS solver is compared against Newton's method on benchmark problems are given.
In particular, for a refined version of the SAIGUP model, the FAS solver is at least 35\% faster than Newton's method for time steps with a CFL number greater than 10.

\end{abstract}


\end{frontmatter}



\section{Introduction} \label{sec:intro}

Numerical simulation is an essential tool to better understand and manage subsurface systems in a wide range of applications including CO$_2$ storage, geothermal energy production, and underground hydrogen storage. 
Solving the partial differential equations (PDEs) governing these systems is challenging due to the nonlinearity of subsurface transport processes and the high heterogeneity of geological porous media.
To avoid unpractical time step size restrictions in the presence of large Courant-Friedrichs-Lewy (CFL) numbers, the temporal discretization of choice is often the unconditionally stable fully implicit method (FIM). 
This approach is computationally expensive as it requires solving large, non-symmetric, ill-conditioned linear systems at each nonlinear iteration. 
Fast nonlinear solvers exhibiting robust convergence properties for large time steps are therefore desirable to reduce the computational burden of the simulations.

Newton-Krylov methods are widely used to solve the nonlinear systems arising from the FIM.
This approach is often combined with linear multigrid \cite{brandt77} to accelerate the Krylov algorithm employed at each Newton iteration.
However, the performance of Newton's method can be undermined by slow nonlinear convergence, especially for large time steps and/or poor initial guesses.
To address this limitation, significant efforts have been invested to enlarge the Newton convergence radius using damping heuristics \cite{jenny2009, wang2013, li2015,moyner2017}, smoother finite volume discretizations \cite{lee2015hybrid,hamon2016implicit,moncorge2020consistent,bosma2022smooth}, and reordering methods \cite{kwok2007, natvig2008, hamon2016ordering, kelmetstal2020reordering}. 
Nonlinear preconditioners based on additive and multiplicative Schwarz preconditioned inexact Newton (A/MSPIN) \cite{cai2002, liu2015, dolean2016} have also been proposed to accelerate nonlinear convergence of multiphase flow in porous media \cite{skogestad2013, skogestad2016, kelmetstal2020schwarz, n2023comparison}. 
In this work, we focus on nonlinear multigrid based on a full approximation scheme (FAS) \cite{brandt77}. 
Instead of using multigrid as a preconditioner for the linear systems, our approach leverages the multigrid concept at the nonlinear level to design a scalable nonlinear solver with robust convergence properties for challenging porous media flow applications \cite{christensen16, christensen18, toft18, fas-spectral-diffusion,fas-two-phase}.

The present article aims at extending the applicability of the nonlinear multigrid solver presented in previous work \cite{fas-spectral-diffusion,fas-two-phase} to the presence of complex line source/sink boundary conditions representing multiperforation wells.
We consider here the widely used Peaceman well model \cite{Pea78} to approximate the difference between the wellbore pressure and the pressure in the perforated reservoir cells.
The accurate representation of multiperforation wells in reservoir simulation is a prerequisite to handle real-field problems but remains a challenge in multilevel solution strategies.
These solvers rely on a coarse problem representation in which the complex radial flow dynamics taking place in the near-well region may be poorly captured, resulting in crippled nonlinear convergence.
This limitation, tackled here in the nonlinear multigrid framework, was investigated in previous work for other classes of multilevel methods.
In the context of multiscale finite volume (MSFV), Wolfsteiner \textit{et al.} \cite{wolfsteiner06} modified the definition of the basis functions in the near-well region to better model radial flow at the coarse level.
Multiple improvements have also been proposed to accommodate the presence of wells in the multiscale mixed finite element method (MsMFEM) \cite{arbogast2002two,chen2003numerical,aarnes2004use,aarnes06}.
The findings of Ligaarden \cite{ligaarden-thesis} and Skaflestad and Krogstad \cite{skaflestad2008multiscale} are particularly relevant to this work.
They highlighted the impact of the well location in the perforated coarse cells on the ability of the basis functions to faithfully approximate radial flow near the well.
Specifically, they showed that wells located close the boundaries of coarse cells were inaccurately represented at the coarse level.

We leverage the latter findings in the construction of the FAS nonlinear multigrid solver presented in this work. 
We extend the mixed fractional-flow velocity-pressure-saturation for two-phase flow considered in \cite{fas-two-phase} to include reservoir-well fluxes at the perforations and the well (rate or pressure) control equations.
In the discrete systems, the set of primary unknowns is enlarged to include well perforation total flux and well pressure---in addition to the reservoir total flux, reservoir pressure, and reservoir saturation representing the reservoir state.
The aggregation algorithm is key to obtain basis functions able to represent near-well radial flow at the coarse levels and preserve a robust nonlinear convergence in the presence of wells.
Following the observations of Ligaarden \cite{ligaarden-thesis}, this is achieved by modifying the coarsening algorithm of \cite{fas-two-phase} to ensure that wells are away from the coarse cell boundaries.
We also improve the intergrid operators described in \cite{fas-two-phase} to transfer both well and reservoir data between levels.
Importantly, the modified multigrid solution strategy remains algebraic and applicable to fully unstructured meshes.
A solution strategy for the reservoir-well linear systems is also designed.
We consider three challenging numerical tests that include several multiperforation wells to illustrate the behavior of the coarsening algorithm around wells.
These benchmark cases demonstrate the excellent nonlinear behavior of FAS nonlinear multigrid compared to Newton's method, especially for large time steps.

We review the governing equations of two-phase flow in porous media in Section \ref{sec:model_problem}. 
Section~\ref{sec:discrete_problem} describes the finite-volume discretization and the Peaceman model for multiperforation wells. 
The treatment of the wells in FAS is presented in Section~\ref{sec:full_approximation_scheme} and nonlinear smoothing is discussed in Section~\ref{sec:nonlinear_smoothing}.
Challenging numerical examples illustrating the robustness of the FAS nonlinear multigrid solver are in Section~\ref{sec:numerical_examples}.

\section{Model problem}\label{sec:model_problem}

We consider a two-phase flow and transport problem involving two immiscible and incompressible phases---a wetting phase, $\wetting$, and a non-wetting phase, $\nonWetting$---flowing in an incompressible porous medium.
Throughout the paper, we use subscript $\wetting$ to indicate wetting-phase quantities, and use superscript $\well$ to indicate well-related quantities.
We focus on a mixed fractional-flow velocity-pressure-saturation formulation.
Gravitational effects are ignored and will be considered in a future publication.
We also neglect capillary forces, a frequent assumption in many practical engineering applications.
Therefore, the pressure is the same for both phases, i.e. $p_\wetting = p_\nonWetting = p$.
In this work, using the saturation constraint $\sum_{\alpha=\{\wetting,\nonWetting\}} s_{\alpha} = 1$, we treat the wetting-phase saturation as primary unknown and denote it from now on as $s = s_\wetting$.

For a simply-connected polyhedral domain  $\Omega \in \mathbb{R}^3$ and time interval $\TimeDomain := (T_0, T_f)$, with $T_0$ and $T_f$ the initial and final time, respectively, the strong form of the initial/boundary value problem (IBVP) consists in finding the total Darcy velocity $\tensorOne{v}: \Omega \times \TimeDomain \rightarrow \mathbb{R}^3$, the pressure $p: \Omega \times \TimeDomain \rightarrow \mathbb{R}$, and the wetting-phase saturation $s: \Omega \times \TimeDomain \rightarrow \mathbb{R}$ such that \cite{aziz79}:
%
\begin{subequations}
\begin{align}
	&\frac{1}{\lambda (s)} \mathbb{K}^{-1} \cdot \tensorOne{v} + \nabla p
	= 0
	&& \text{in} \; \Omega \times \TimeDomain
	&& \mbox{(total Darcy velocity)} ,
	\label{eq:total_darcy} \\
	&\nabla \cdot \tensorOne{v}
	=
	q^\well(p, s)
	&& \text{in} \; \Omega \times \TimeDomain
	&& \mbox{(total volume conservation)},
	\label{eq:total_mass_conservation} \\
	&\phi \frac{\partial s}{\partial t} + 
  \nabla\cdot [ f_\wetting (s) \tensorOne{v} ]
	=
	q_\wetting^\well(p,s),
	&& \text{in} \; \Omega \times \TimeDomain
	&& \mbox{(wetting-phase volume conservation)} ,
	\label{eq:wetting_mass_conservation}
\end{align}
\label{eq:model}\null
\end{subequations}
%
where
\begin{itemize}
  \item $\lambda(s) = \sum_{\alpha=\{\wetting,\nonWetting\}} \lambda_{\alpha}(s)$ is the total mobility, with the corresponding phase-based quantities defined as the ratio of relative permeability, $k_{r,\alpha}$, to viscosity, $\mu_\alpha$, i.e. $\lambda_\alpha(s) := k_{r,\alpha}(s) / \mu_\alpha$. Various constitutive relationships for $k_{r,\alpha}$ will be considered in our numerical examples. Note that using standard assumptions on the phase mobilities, the total mobility is bounded away from zero;
  \item $\mathbb{K}$ and $\phi$ are the medium absolute permeability tensor and porosity, respectively;
  \item $q^\well(p,s) = \sum_{\alpha=\{\wetting,\nonWetting\}} q^\well_{\alpha}(p,s)$ is the total volumetric source/sink per unit volume, with $q^\well_{\alpha}(p,s)$ the corresponding phased-based quantity, which are driven by wells and will be defined in more details below.
  \item $f_\wetting(s) := \lambda_\wetting(s) / \lambda(s)$ is the fractional flow function.
\end{itemize}
Without loss of generality, in our simulations the domain boundary, $\partial \Omega$, is always subject to no-flow boundary conditions.
This represents a natural assumption when simulating closed-flow systems, such as depleted reservoirs in which CO$_2$ is injected for permanent storage.
To ensure uniqueness of the pressure solution, we prescribe a datum value for pressure internally in the domain through sink and/or source terms.
The formulation is completed by appropriate initial conditions for $\tensorOne{v}$, $p$, and $s$.

\section{The discrete problem}\label{sec:discrete_problem}

The system of PDEs \eqref{eq:model} is discretized by a cell-centered two-point flux approximation (TPFA) finite-volume (FV) method \cite{EymGalHer00} on a conforming triangulation of the domain, combined with the backward Euler (fully implicit) time-stepping scheme.

\subsection{Finite Volume discretization}\label{sec:FV_tpfa}

First, we introduce some notation.
Let $\mathcal{T}$ be the set of cells in the computational mesh such that $\overline{\Omega} = \sum_{\tau \in \mathcal{T}} \overline{\tau}$.
For a cell $\tau_K \in \mathcal{T}$, with $K$ a global index, let $|\tau_K|$ denote the volume, $\partial \tau_{K} = \overline{\tau}_{K} \setminus \tau_{K}$ the boundary, $\tensorOne{x}_{K}$ the barycenter, and $\tensorOne{n}_{K}$ the outer unit normal vector associated with $\tau_K$.
Also, let $\mathcal{T}_\well$ be the subset of $\mathcal{T}$ containing cells that are perforated by some wells.
Let $\mathcal{E}$ be the set of internal faces in the computational mesh included in $\Omega$.
An internal face $\varepsilon$ shared by cells $\tau_K$ and $\tau_L$ is denoted as $\varepsilon_{K,L} = \partial \tau_K \cap \partial \tau_L$, with the indices $K$ and $L$ such that $K < L$.
%
%
The area of a face is $|\varepsilon|$.
A unit vector $\tensorOne{n}_{\varepsilon}$ is introduced to define a unique orientation for every face, and we set $\tensorOne{n}_{\varepsilon} = \tensorOne{n}_{K}$.
To indicate the mean value of a quantity $(\cdot)$ over a face $\varepsilon$ or a cell $\tau_K$, we use the notation $(\cdot)_{\left| \right. \varepsilon}$ and $(\cdot)_{\left| \right. K}$, respectively.
Let $T_0 = t_0 < t_1 < \cdots < t_n = T_f$ be a partition of the time domain $\TimeDomain$.
The discrete (finite-difference) approximation to a time-dependent quantity $\chi(t_m)$ at time $t_m$ is denoted by $\chi^m$.
Also, we define the time step size $\Delta t_m := t_m - t_{m-1}$.

We consider a piecewise-constant approximation for both pressure and saturation.
For each cell $\tau_{K} \in \mathcal{T}$, we introduce one pressure, $p_K$, and one saturation, $s_K$, degree of freedom, respectively.
We denote by $\sigma_{\varepsilon}$ the numerical flux approximating the total Darcy flux through an internal face $\varepsilon=\varepsilon_{K,L}$, i.e. $\sigma_{\varepsilon} \approx \int_{\varepsilon_{K,L}} \tensorOne{v} \cdot \tensorOne{n}_{\varepsilon} \mathrm{d}\Gamma$, such that:
%
\begin{equation}
  \left( \frac{1}{\lambda(s_K) \overline{\Upsilon}_{K,\varepsilon}} + \frac{1}{\lambda(s_L) \overline{\Upsilon}_{L,\varepsilon}} \right) \sigma_{\varepsilon} -  (p_{K} - p_{L}) = 0,  
  \label{eq:TPFA_flux}
\end{equation}
%
where $\overline{\Upsilon}_{K,\varepsilon} (s_K)$ and $\overline{\Upsilon}_{L,\varepsilon} (s_L)$ are the constant (geometric) one-sided transmissibility coefficients, defined as \cite{lie19}
%
\begin{align}
  \overline{\Upsilon}_{i,\varepsilon} &=
  |\varepsilon_{K,L}| \frac{\tensorOne{n}_{i} \cdot \mathbb{K}_{\left| \right. i} \cdot (\tensorOne{x}_{\varepsilon} - \tensorOne{x}_{i})}{||\tensorOne{x}_{\varepsilon} - \tensorOne{x}_{i}||_2^2},
  &
  i &= \{K, L \},
  \label{eq:half-transmissibility}
\end{align}
%
with $\tensorOne{x}_{\varepsilon}$ a collocation point introduced for every $\varepsilon \in \mathcal{E}$ to enforce point-wise pressure continuity across interfaces.

\noindent
The approximation of the wetting-phase Darcy flux through $\varepsilon = \varepsilon_{K,L} \in \mathcal{E}$ in the discrete form of Eq.~\eqref{eq:wetting_mass_conservation} relies on using single-point upstream weighting (SPU) according to the sign of $\sigma_{\varepsilon}$, namely
%
\begin{align}
  f_\wetting^{\text{upw}}(s_K,s_L) \sigma_{\varepsilon} 
  &\approx \int_{\varepsilon_{K,L}} f_\wetting(s)\tensorOne{v}\cdot\tensorOne{n}_{\varepsilon}\mathrm{d}\Gamma,
  &
  f_\wetting^{\text{upw}}(s_K,s_L)
  &=
  \begin{cases}
    f_\wetting(s_K), & \text{if } \sigma_{\varepsilon} > 0, \\
    f_\wetting(s_L), & \text{otherwise}.
  \end{cases}
  \label{eq:upwinding}
\end{align}
%

\subsection{A standard well model with multiple perforations}

\begin{figure}
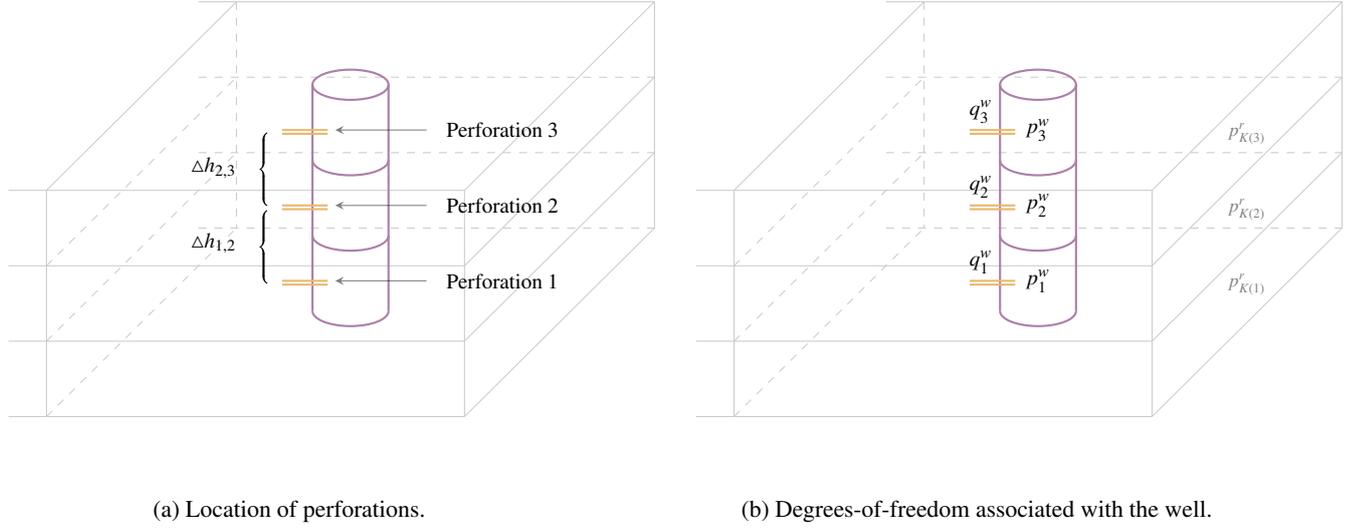

\centering
  \begin{subfigure}[b]{0.45\textwidth}
    \centering
    \include{./pics/well}
    \caption{Location of perforations.}
    \label{fig:well_example}
  \end{subfigure}
  \hfill
  \begin{subfigure}[b]{0.45\textwidth}
    \centering
    \include{./pics/well_dofs}
    \caption{Degrees-of-freedom associated with the well.}
    \label{fig:well_dof_example}
  \end{subfigure}
  \caption{An example of a 3-perforation well.}
  \label{fig:multi-perforation-well}
\end{figure}

Source and sink terms in Eqs.~\eqref{eq:total_mass_conservation}-\eqref{eq:wetting_mass_conservation} are used to simulate the effect of injection and production wells.
We employ a conventional Peaceman well model \cite{Pea78}, which relates well control parameters, such as bottomhole pressure (BHP), to flow rates through the wellbore \cite{Pea78}.
We assume each well segment to be vertical, with a single perforation connected to the centroid of a reservoir cell.
Also, without loss of generality, we restrict ourselves to rate-controlled injection wells and BHP-controlled production wells.
Consider a well with $n_{perf}$ perforations, see Figure~\ref{fig:multi-perforation-well} for an example.
The pressure values $p^\well_i$ at the perforations $i$ are the unknowns.
We can write $n_{perf}$ equations for the well, where the first $n_{perf} - 1$ equations describe the pressure relation between
two neighboring perforations. Specifically, the expression can be written as
\begin{equation}
p^\well_{i+1} - p^\well_i = \frac{1}{2} (\rho^\well_{i} + \rho^\well_{i+1} ) g \Delta h_{i,i+1} \qquad i = 1, \dots, n_{perf} - 1,
\label{eq:well_pressure}
\end{equation}
where $\rho^\well_{i}$ stands for the density of the fluid mixture in the wellbore around perforation $i$, $g$
denotes the gravity constant, and $\Delta h_{i,i+1}$ is the height difference between perforations $i$ and $i + 1$.
In the current paper, we ignore the effect of gravity (i.e., $g = 0$), so \eqref{eq:well_pressure} implies that the perforation pressures for the same well are the same. Hence, we have effectively only one pressure unknown associated with a well $\well$, and it is denoted by
\[
p^\well := p^\well_1 = p^\well_2 = \cdots = p^\well_{n_{perf}}.
\]
%
The phase volumetric rate $q^\well_{\alpha, i}$ through perforation $i$ is calculated with the Peaceman model:
\begin{equation}
  q_{\alpha, i}^\well = \lambda_{\alpha}(s_{\alpha, K(i)}) WI_i \left(p^{\res}_{K(i)} - p^\well_i \right) = WI_i\; \lambda_{\alpha}(s_{\alpha, K(i)}) \left(p^{\res}_{K(i)} - p^\well \right) \qquad i = 1, \dots, n_{perf}.
  \label{eq:phase_volumetric_rate}
\end{equation}
Here, $K(i)$ is the index of the reservoir cell in $\mathcal{T}$ where perforation $i$ is located,
$WI_i$ is the well index describing the transmissibility between the wellbore and the perforated reservoir cell, 
and $p^{r}_{K(i)}$ is the pressure in the perforated cell.
A comprehensive presentation of well models and well index calculation can be found in \cite{CheHuaMa06}.
The source/sink term in \eqref{eq:model} is then defined as
\begin{equation}
  q_{\alpha}^\well(p, s) := \sum_\well\sum_{i = 1}^{n_{perf}} q_{\alpha,i}^\well \delta ( \tensorOne{x} - \tensorOne{x}_{K(i)} )
  \label{eq:source_def}
\end{equation}
where $\delta ( \cdot )$ is the Dirac function and $\tensorOne{x}_{K(i)}$ is the barycenter of the perforated cell $\tau_{K(i)}$.
We introduce a new unknown $\sigma_{i}^\well$ for each well cell $\tau_{K(i)}\in\mathcal{T}_\well$:
\begin{equation}
  \sigma_{i}^\well 
  := \int_{\tau_{K(i)}} -q^\well(p, s) = \int_{\tau_{K(i)}} \sum_{\alpha=\{\wetting,\nonWetting\}} -q^\well_{\alpha}(p,s) \qquad i = 1, \dots, n_{perf}.
  \label{eq:well_flux_def}
\end{equation}
Hence, by \eqref{eq:phase_volumetric_rate}, \eqref{eq:source_def}, and \eqref{eq:well_flux_def}
\begin{equation}
  \frac{1}{\lambda_{T}(s_{K(i)})WI_i}\sigma_{i}^\well = \frac{1}{\lambda_{T}(s_{K(i)})WI_i} \left( -q_{\wetting, i}^\well - q_{\nonWetting, i}^\well \right) = \left( p^\well - p^\res_{K(i)} \right) \qquad i = 1, \dots, n_{perf}
  \label{eq:total_volumetric_rate}
\end{equation}
and
\begin{equation}
  f_\wetting(s_{K(i)})\sigma_{i}^\well = q_{\wetting, i}^\well \qquad i = 1, \dots, n_{perf}.
\end{equation}

\subsection{Control equations}

In addition to the $n - 1$ pressure relation equations \eqref{eq:well_pressure}, there is one
more equation for each well, which is the well constraint equation.
This constraint equation reflects the physical control strategy on a well.
The most common constraint equations are phase rate control and wellbore pressure control.
For a well with BHP control, the constraint equation can be written as follows:
\begin{equation}
p^{\well} - p^{\well,target} = 0,
\end{equation}
where 
$p^{w,target}$ is the user-specified operating pressure.
For a well with phase rate control, the constraint equation for a well is:
\begin{equation}
\sum^{n_{perf}}_{i = 1} q_{\alpha, i}^\well -  q^{\well,target}_\alpha = 0,
\end{equation}
where $q^{\well,target}_\alpha$ is the specified injection rate of phase $\alpha$.

%
%
%

\subsection{Structure of the discrete problem}
Let $\mathcal{W}$ be the set of wells. Introducing coefficient vectors
$\Vec{\sigma}^{\res} = (\sigma_{\varepsilon})_{\varepsilon \in \mathcal{E}}$,
$\Vec{\sigma}^{\well} = (\sigma_{K}^{\well})_{\tau_K \in \mathcal{T}_\well}$,
$\Vec{p}^{\res} = (p_{K})_{\tau_{K} \in \mathcal{T}}$,
$\Vec{p}^{\well} = (p^{\well}_w)_{w \in \mathcal{W}}$
and
$\Vec{s} = (s_{K})_{\tau_{K} \in \mathcal{T}}$
that contain the unknown degrees of freedom at time $t = t_m$ (i.e. face fluxes, perforation fluxes, cell pressures, well pressures, and cell saturations) the algebraic form associated with the IBVP \eqref{eq:model} can be stated as follows: 

\begin{equation}
\Vec{r}(\Vec{x}) 
:=
\begin{bmatrix}
\Vec{r}_{\sigma^\res}(\Vec{x}) \\
\Vec{r}_{\sigma^\well}(\Vec{x}) \\
\Vec{r}_{p^\res}(\Vec{x}) \\
\Vec{r}_{p^\well}(\Vec{x}) \\
\Vec{r}_s(\Vec{x})  
\end{bmatrix}
: =
\begin{bmatrix}
M^\res(\Vec{s})\Vec{\sigma}^{\res} - (D^{\res, \res})^T \Vec{p}^{\res} \\
M^\well(\Vec{s})\Vec{\sigma}^{\well} - (D^{\res, \well})^T \Vec{p}^{\res} - (D^{\well, \well})^T \Vec{p}^{\well} - \Vec{g}^{\well}\\
D^{\res, \res}\Vec{\sigma}^{\res} + D^{\res, \well}\Vec{\sigma}^{\well} \\
D^{\well, \well}\Vec{\sigma}^{\well} - \Vec{f}^{\well} \\
T(\Vec{\sigma}^{\res}, \Vec{\sigma}^{\well}, \Vec{s}) - (\Delta t_m)^{-1}W\Vec{s}^{m-1}
\end{bmatrix} = \Vec{0},
\quad \text{where} \quad
\Vec{x}
:=
\begin{bmatrix}
\Vec{\sigma}^{\res} \\
\Vec{\sigma}^{\well} \\
\Vec{p}^{\res} \\
\Vec{p}^{\well} \\
\Vec{s}
\end{bmatrix},
\label{eq:discrete_problem}
\end{equation}
%
$(D^{\res, \res})^T$ computes for the pressure difference between adjacent cells as in \eqref{eq:TPFA_flux}, $[(D^{\res, \well})^T \; (D^{\well, \well})^T]$ computes for the pressure difference between wells and well cells as in \eqref{eq:phase_volumetric_rate}, $M^\res(\Vec{s})$ has the coefficients of $\sigma_\varepsilon$ in \eqref{eq:TPFA_flux} on the diagonal, $M^\well(\Vec{s})$ has the inverse of the scaled well indices in \eqref{eq:TPFA_flux} on the diagonal, $T$ comes from the discretization of \eqref{eq:wetting_mass_conservation} as in \cite{fas-two-phase},
$\Vec{g}^{\well}$ contains the BHP, $\Vec{f}^{\well}$ contains the injection rates, $\Vec{s}^{m-1}$ is the coefficient vector of the saturation solution at the previous time step , and
\begin{equation}
T(\Vec{\sigma}^{\res}, \Vec{\sigma}^{\well}, \Vec{s}) = (\Delta t_m)^{-1}W \Vec{s} + D^{\res,\res} \text{diag}\left( \Vec{\sigma}^{\res} \right) U(\Vec{\sigma}^{\res}) f_w(\Vec{s}) + D^{\res,\well} \text{diag}\left( \Vec{\sigma}^{\well} \right) U(\Vec{\sigma}^{\well}) f_w(\Vec{s}^{\well}).
\label{eq:discrete_transport}
\end{equation}
Here, $\text{diag}\left(\Vec{\sigma}\right)$ is the diagonal matrix created from the entries of the argument vector, and $U(\Vec{\sigma}^\res)$ is the upwind operator selecting for each connection the appropriate upstream value from the input vector $f_w(\Vec{s})$, which contains the fractional flow function values evaluated in each cell, with $\Vec{s}^{\well}$ the saturation value in the wells.

\begin{figure}[ht]
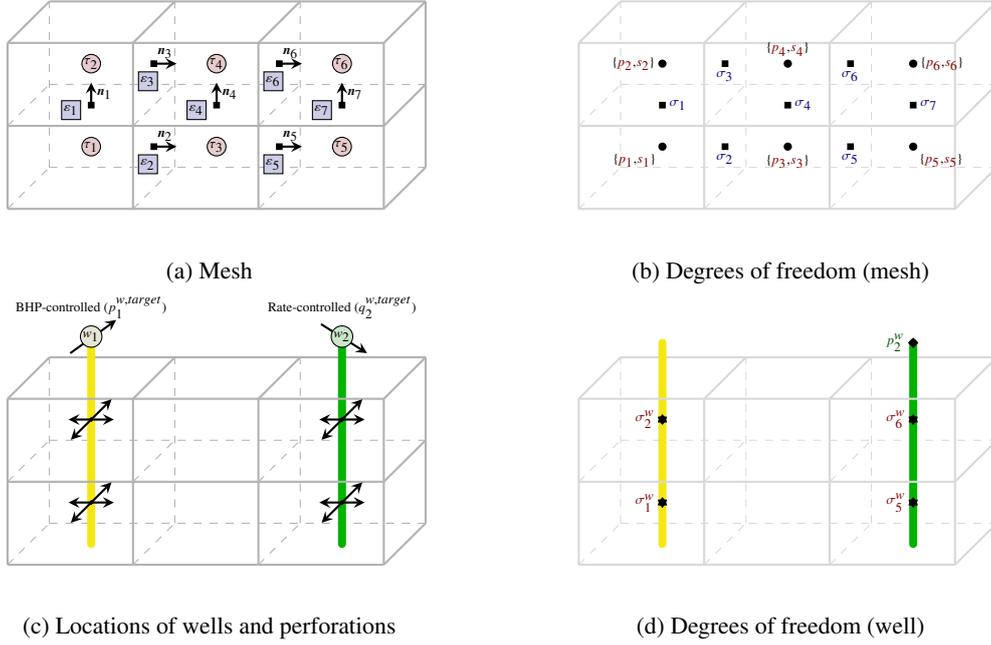

  \tiny
  \centering
  \begin{subfigure}[b]{0.33\textwidth}
    \centering
    \include{pics/app_mesh_geometry}
    \caption{Mesh}
    \label{fig:mesh_sketch}
  \end{subfigure}
  \hspace{20mm}
  \begin{subfigure}[b]{0.33\textwidth}
    \centering
    \include{pics/app_mesh_dofs}    
    \caption{Degrees of freedom (mesh)}
    \label{fig:mesh_dofs}
  \end{subfigure}
  \\
  \begin{subfigure}[b]{0.33\textwidth}
    \centering
    \include{pics/app_mesh_wells}    
    \caption{Locations of wells and perforations}
    \label{fig:mesh_wells}
  \end{subfigure}
  \hspace{20mm}
  \begin{subfigure}[b]{0.33\textwidth}
    \centering
    \include{pics/app_well_dofs}    
    \caption{Degrees of freedom (well)}
    \label{fig:well_dofs}
  \end{subfigure}
 \caption{Sketch of a well-driven flow using a mesh consisting of six cells.  The domain boundary is subject to no-flow conditions everywhere.  The locations of the rate-controlled injection well (green) and the BHP-controlled production well (yellow) are shown in (c). Note that the numbering of $\sigma^{\well}_i$ is based on the definition in \eqref{eq:well_flux_def}.}
 \label{fig:well_driven_flow}
\end{figure}

\begin{figure}[ht]
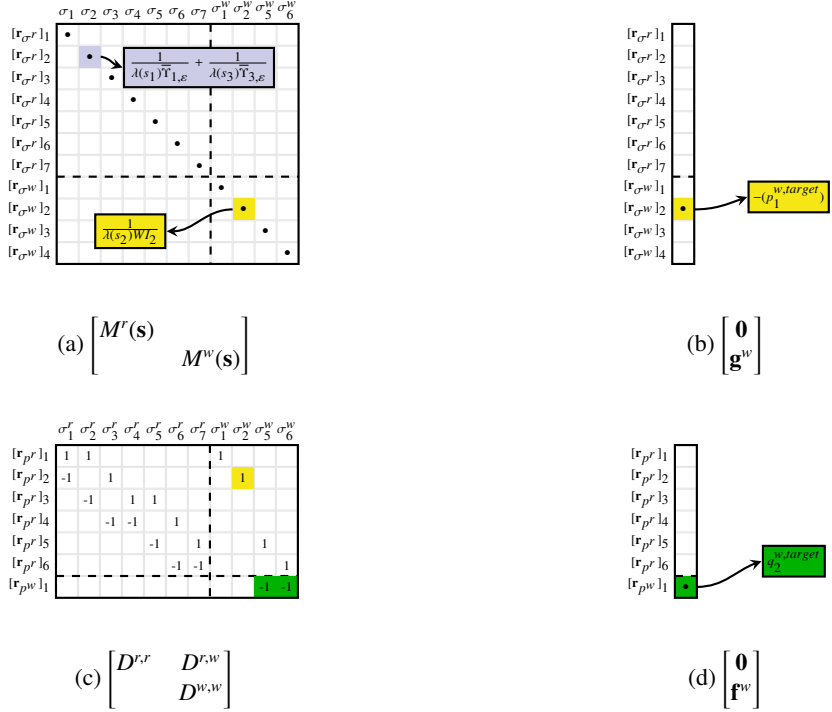

  \tiny
  \centering
  \begin{subfigure}[b]{0.33\textwidth}
    \centering
    \include{pics/app_mat_M}
    \caption{$\begin{bmatrix} M^\res(\Vec{s}) \\ & M^\well(\Vec{s}) \end{bmatrix}$}
    \label{fig:example_M}
  \end{subfigure}
  \hspace{20mm}
  \begin{subfigure}[b]{0.33\textwidth}
    \centering
    \include{pics/app_vec_g}
    \caption{$\begin{bmatrix} \Vec{0} \\ \Vec{g}^\well \end{bmatrix}$}
    \label{fig:example_g}
  \end{subfigure}
  \\
  \vspace{5mm}
  \begin{subfigure}[b]{0.33\textwidth}
    \centering
    \include{pics/app_mat_D}
    \caption{$\begin{bmatrix} D^{\res,\res} & D^{\res,\well} \\ & D^{\well,\well}  \end{bmatrix}$}
    \label{fig:example_D}
  \end{subfigure}
  \hspace{20mm}
  \begin{subfigure}[b]{0.33\textwidth}
    \centering
    \include{pics/app_vec_f}
    \caption{$\begin{bmatrix} \Vec{0} \\ \Vec{f}^\well \end{bmatrix}$}
    \label{fig:example_f}
  \end{subfigure}
 \caption{Components of the discrete problem for Fig.~\ref{fig:well_driven_flow}. The terms in yellow are involved in the discrete pressure constraint for well 1, while the terms in green correspond to the discrete rate constraint for well 2.}
 \label{fig:discrete_components_example}
\end{figure}

To illustrate the components in \eqref{eq:discrete_problem}, we consider a simple well-driven flow example with the mesh and degrees-of-freedom depicted in Fig.~\ref{fig:well_driven_flow}.
The corresponding discrete operators and vectors in \eqref{eq:discrete_problem} for Fig.~\ref{fig:well_driven_flow} are shown in Fig.~\ref{fig:discrete_components_example}.

\section{Full approximation scheme}\label{sec:full_approximation_scheme}

In this section, we propose a nonlinear multigrid solver for the discrete nonlinear system \eqref{eq:discrete_problem} based on the Full Approximation Scheme (FAS) \cite{brandt77,henson2003multigrid}.
We start by giving a high-level overview of FAS and its essential components.
First, we will need three intergrid transfer operators---namely, an interpolation operator $\blkInterpolate$, a restriction operator $\blkRestrict$, and a projection operator $\blkInject$.
In particular, $\blkInterpolate$ and $\blkInject$ satisfy
\begin{equation}
\blkInject\blkInterpolate = \blkMat{I}^{\ell+1},
\end{equation}
where $\blkMat{I}^{\ell+1}$ is the identity operator on the level $\ell+1$.
\begin{rmk}[Abuse of terminology]

$\blkInject$ is not a projection according to the usual definition of projections. Nevertheless, following the discussion in \cite[Remark~6]{fas-spectral-diffusion}, $\blkInject$ will be referred to as a projection with an abuse of terminology.

\end{rmk}

We use the convention that level $\ell = 0$ refers to the finest level (i.e., the original problem), and a larger value of $\ell$ means a coarser level.
Moreover, a hierarchy of nonlinear operators $\left\{ \Vec{r}^{\ell}(\Vec{x}^\ell) \right\}_{\ell = 0}^{\mathcal{L}-1}$ approximating $\Vec{r}(\Vec{x})$ will need to be constructed.
Lastly, the approximated solution is updated at each level based on some smoothing step, denoted ``{\tt NonlinearSmoothing}."
A typical step at level $\ell$ in the full approximation scheme multigrid is stated in Algorithm~\ref{alg:fas}, where $n_s^\ell$ is the number of smoothing steps at level $\ell$.
The backtracking procedure is described in \cite[Algorithm~1]{fas-spectral-diffusion}.

At time step $m$, the multigrid solver for \eqref{eq:discrete_problem} starts with an initial guess, $\blkVec{x}^{0} := \blkVec{x}^{m-1}$, chosen to be the converged state at the previous time step $m-1$. Then, the solver performs a sequence of nonlinear iterations denoted by the superscript $k$, as follows:
\begin{equation}
\blkVec{x}^{k} =  {\tt NonlinearMG}(0,\, \blkVec{x}^{k-1},\, \Vec{0}), \quad \forall\, k \ge 1,
\end{equation}
until a certain stopping criterion is satisfied.
In the rest of this section, the details of all multigrid cycle components will be discussed.
\begin{algorithm}[t]
\caption{Nonlinear step at level $\ell$ in the Full Approximation Scheme}
\label{alg:fas}
\begin{algorithmic}[1]
\Function{\tt NonlinearMG}{$\ell,\, \blkVec{x}^\ell,\, \blkVec{b}^\ell$}
	\If{$\ell$ is the coarsest level}
		\State $\blkVec{x}^\ell \leftarrow$ \Call{\tt NonlinearSmoothing}{$\ell,\, \blkVec{x}^\ell,\, \blkVec{b}^\ell,\, n_s^\ell$}
	\Else
		\State $\blkVec{x}^\ell \leftarrow$ \Call{\tt NonlinearSmoothing}{$\ell,\, \blkVec{x}^\ell,\, \blkVec{b}^\ell,\, n_s^\ell$} \label{alg:line:nonlinPreSmoothing}
		\State $\blkVec{x}^{\ell+1} \leftarrow \blkInject \blkVec{x}^\ell$
		\State $\blkVec{b}^{\ell+1} \leftarrow \blkVec{r}^{\ell+1}(\blkVec{x}^{\ell+1}) - \blkRestrict(\blkVec{r}^{\ell}(\blkVec{x}^\ell) - \blkVec{b}^\ell )$ 
		\State $\blkVec{y}^{\ell+1} \leftarrow$ \Call{\tt NonlinearMG}{$\ell+1,\, \blkVec{x}^{\ell+1},\, \blkVec{b}^{\ell+1}$}
		\State $\blkVec{x}^\ell \leftarrow$ \Call{\tt Backtracking}{$\blkVec{x}^\ell,\, \blkInterpolate( \blkVec{y}^{\ell+1} - \blkVec{x}^{\ell+1}),\, \theta$} 
		\State $\blkVec{x}^\ell \leftarrow$ \Call{\tt NonlinearSmoothing}{$\ell,\, \blkVec{x}^\ell,\, \blkVec{b}^\ell,\, n_s^\ell$} \label{alg:line:nonlinPostSmoothing}
  	\EndIf
	\State \Return{$\blkVec{x}^\ell$}
 \EndFunction
 \end{algorithmic}
\end{algorithm}

\subsection{Intergrid transfer operators} \label{sec:intergrid_operators}

The interpolation operator $\blkInterpolate$ and the projection operator $\blkInject$ are block-diagonal, composed of the corresponding operators for the flux at faces, the perforation flux, the cell-center pressure, the well pressure, and the saturation unknowns:
\begin{equation}
\blkInterpolate = \begin{bmatrix}
P_{\sigma^\res} \\ & P_{\sigma^\well} \\ & & P_{p^\res} \\ & & & P_{p^\well} \\ & & & & P_s
\end{bmatrix}.
\label{eq:prolongation}
\end{equation}
The restriction operator $\blkRestrict$ and injection operator $\blkInject$ have a similar structure and notation.
The restriction operator $\blkRestrict$ is taken as the transpose of the interpolation operator $\blkInterpolate$: 
\[
\blkRestrict := \blkInterpolate^T.
\]

To define our interpolation operators, we first form a nested hierarchy of grids $\{\grid^\ell\}_{\ell=0}^{\mathcal{L}-1}$ by aggregating fine grid cells in $\grid^0 := \grid$.
Starting with $\ell = 0$, we consider the cell-connectivity graph of $\grid^\ell$, where each cell (respectively face) in $\grid^\ell$ is a vertex (respectively edge) of the graph.
Based on the cell-connectivity graph, we add wells as additional vertices, and well perforations (connections between a well and a reservoir cell) as additional edges to the graph, and we call this the cell-well-connectivity graph.
Based on the cell-well-connectivity graph, contiguous aggregates of cells are formed by using a graph partitioner (e.g., METIS \cite{karypis1998fast}).
Since the wells require a specific treatment, we do not aggregate them with the other cells. 
To this end, we identify them in the cell-well-connectivity and make each well an individual aggregate by itself.
These aggregates are the ``cells" (which have irregular shapes) in the coarser-level grid $\grid^{\ell+1}$.
A coarser-level face is also naturally formed by collecting the fine faces sharing a pair of adjacent aggregates.
The set of faces on level $\ell$ is denoted by $\faces^\ell$.
This process is repeated until the coarsest grid $\grid^{\mathcal{L}-1}$ is formed.

For the construction of the components of $\blkInterpolate$ and $\blkInject$, we use the lowest order graph-based multilevel coarsening method in \cite{ml-spectral-coarsening} applied to the cell-well-connectivity graph (see also \cite{fas-two-phase}).

\subsection{The nonlinear problem on each level} \label{sec:coarse_operators}

On the fine level $\ell = 0$, let $\Vec{r}^{0}(\Vec{x}^0) := \Vec{r}(\Vec{x}^0)$.
With the interpolation operator $\blkInterpolate$ and restriction operator $\blkRestrict$, the nonlinear operators on coarse levels are defined recursively as
\begin{equation}
\Vec{r}^{\ell+1}(\Vec{x}^{\ell+1}) := \blkRestrict\Vec{r}^{\ell}\left( \blkInterpolate\Vec{x}^{\ell+1} \right)
\label{eq:level_operator}
\end{equation}
where
\begin{equation}
T^{\ell+1}(\Vec{\sigma}^{\res, \ell+1}, \Vec{\sigma}^{\well, \ell+1}, \Vec{s}^{\ell+1}) := \restrict{s} T^{\ell}\left( \interpolate[\res]{\sigma} \Vec{\sigma}^{\res, \ell+1}, \interpolate[\well]{\sigma} \Vec{\sigma}^{\well, \ell+1},  \interpolate{s} \Vec{s}^{\ell+1} \right), 
\end{equation}
The nonlinear problem on level $\ell$ is 
\begin{equation}
\Vec{r}^{\ell}(\Vec{x}^{\ell}) - \Vec{b}^{\ell} = \Vec{0}^{\ell}
\label{eq:level_problem}
\end{equation}
where $\Vec{b}^{\ell}$ is defined recursively as in Algorithm~\ref{alg:fas} with $\Vec{b}^{0} = \Vec{0}^0$, a vector of all zero on the finest level.
Note that \eqref{eq:level_operator} is a conceptual definition of the coarse operators. 
In practice, due to scalability concerns, we do not want the evaluation of $\Vec{r}^{\ell+1}(\Vec{x}^{\ell+1})$ during the multigrid cycle to involve computations on the finer level $\ell$.
To achieve this, we have to construct and store some coarse operators and vectors during the setup phase of the multigrid solver.
For $D^{\ell+1}, \Vec{g}^{\ell+1}, \Vec{f}^{\ell+1}$, and $\Vec{h}^{\ell+1}$, the construction is straightforward.
The main issue is in the evaluation of the nonlinear components $M^{\ell+1}\left( \Vec{s}^{\ell+1} \right)$ and $T^{\ell+1}(\Vec{\sigma}^{\res, \ell+1}, \Vec{\sigma}^{\well, \ell+1}, \Vec{s}^{\ell+1})$, where the details are discussed in \cite{fas-two-phase}.

\subsection{Aggregation around wells}
\label{sec:aggregation_around_wells}

\begin{figure}[h]
  \begin{subfigure}[b]{0.45\textwidth}
    \centering
    \includegraphics[scale=.18]{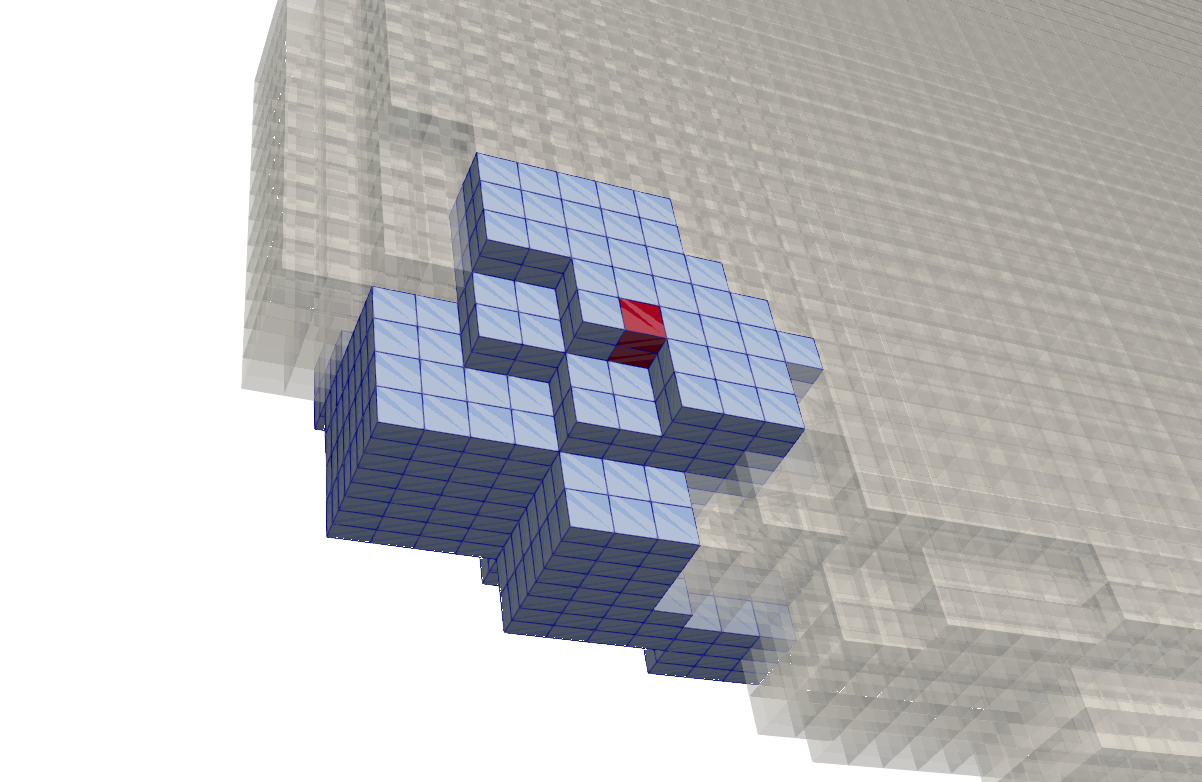}
    \caption{view from top}
  \end{subfigure}
  \hspace{5mm}
  \begin{subfigure}[b]{0.45\textwidth}
    \centering
    \includegraphics[scale=.18]{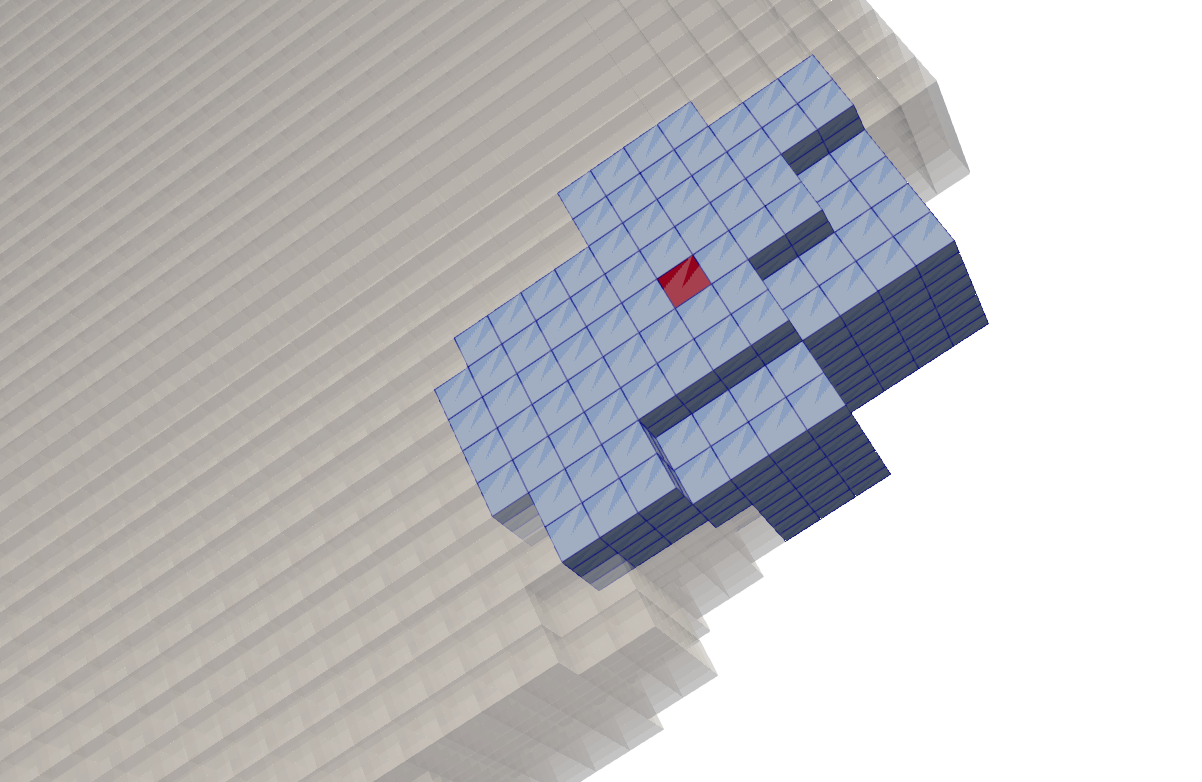}
    \caption{view from bottom}
  \end{subfigure}
\caption{An example aggregate containing a well using Algorithm~\ref{alg:well-part} for the refined Egg model. The fine-level reservoir cells perforated by the well are colored in red.}
\label{fig:well-aggregate}
\end{figure}

It is known that the location of the wells in the coarse elements can greatly affect the quality of the multiscale basis functions associated with the wells. 
For example, when a well is located at a corner or a boundary face of a coarse element, the subgrid feature of the multiscale basis function caused by the well cannot be extended to the neighboring coarse elements, leading to a mismatch with the global solution \cite{ligaarden-thesis}. 
To tackle this issue, we assign heavier weights for edges corresponding to faces nearby the well (a few layers around the well cells) in the cell-connectivity graph, and then feed the modified weighted graph to METIS.
After the initial aggregation is formed, we merge the aggregates that are connected to some well cells of the same well.
This second step ensures that well cells are not in contact with the boundary of the aggregate (unless some face of the well cell is originally on the boundary of the computational domain). 
The overall aggregation algorithm is summarized in Algorithm~\ref{alg:well-part}.
Figure~\ref{fig:graph-weight-example} illustrates the location of edges whose weights are modified in a small graph example when {\tt n\_lay} (number of layers) is 2 in Algorithm~\ref{alg:well-part}.
In our numerical experiments, we set the parameters ``number of layers" to 4 and ``edge weight multiplier" to $10^6$ in Algorithm~\ref{alg:well-part}.
With this approach, we observe that well cells 
stay away from aggregate-to-aggregate interface in the resulting partition. See Figure~\ref{fig:well-aggregate} for an example of the resulting aggregate containing the well cells.

\begin{algorithm}[t]
\caption{Well-location informed partition}
\label{alg:well-part}
\begin{algorithmic}[1]
    \State {\bf Input:} cell-connectivity graph ({\tt conn}), well cell indices ({\tt well\_cells})
    \State {\bf Output:} cell partitioning array ({\tt part})
    \State {\bf Parameters:} number of layers ({\tt n\_lay}), edge weight multiplier ({\tt scale})
    
    {\textcolor{white}{}} 
    
    \State {\tt cell\_layers} $\gets \Pi^{{\tt n\_lay}}_{i=1}${\tt conn} \Comment{find out cells that are {\tt n\_lay} (cell) layers away from well cells}

    \State  {\tt well\_cell\_neighbors} $\gets$ cells connected to {\tt well\_cells} in {\tt cell\_layers}
    \For{cell-cell connection in {\tt conn} (an edge in the sparse matrix graph)}
        \If{both cells of the connection are in {\tt well\_cell\_neighbors}}
            \State  {\tt edge\_weight} $\gets$ {\tt edge\_weight * scale} \Comment{Modify edge weights of {\tt conn} near well}
        \EndIf
    \EndFor

    \State {\tt part\_0} $\gets$ METIS(modified {\tt conn}) \Comment{Produce a partitioning array by calling  METIS}
    \State {\tt part} $\gets$ merge aggregates in {\tt part\_0} that are connected to some well cells of the same well 
    \State \Return{\tt part}
\end{algorithmic}
\end{algorithm}

\begin{figure}
\centering
    \include{./pics/well-aggregation}
    \caption{Location of edges with modified weights in an example cell-connectivity graph if {\tt n\_lay} = 2 in Algorithm~\ref{alg:well-part}.}
    \label{fig:graph-weight-example}
\end{figure}

\section{Nonlinear smoothing}\label{sec:nonlinear_smoothing}

In our FAS scheme, the nonlinear smoothing step on level $\ell$ in Algorithm~\ref{alg:fas} consists of applying a Newton iteration to solve \eqref{eq:level_problem}. The number of Newton iterations is a user-specified parameter. 
As shown in \cite{fas-two-phase}, the coarse-level problem has the same structure as the original fine-level problem \eqref{eq:discrete_problem}. 
%
%
Therefore, the Jacobian on any level has the form:
\begin{align}
    \begin{bmatrix}
		\Mat{A}_{\sigma^{\res} \sigma^{\res}} & \Mat{A}_{\sigma^\res \sigma^\well} & \Mat{A}_{\sigma^\res p^r} &  & \Mat{A}_{\sigma^\res s} \\
		\Mat{A}_{\sigma^\well\sigma^\res} & \Mat{A}_{\sigma^\well\sigma^\well} & \Mat{A}_{\sigma^\well p^\res} & \Mat{A}_{\sigma^\well p^\well} & \Mat{A}_{\sigma^\well s} \\
		\Mat{A}_{p^\res\sigma^\res} & \Mat{A}_{p^r\sigma^\well} &  &  &  \\
		 & \Mat{A}_{p^\well\sigma^\well} &  &  &  \\
		\Mat{A}_{s\sigma^\res} & \Mat{A}_{s\sigma^\well} &  &  & \Mat{A}_{ss} \\
	\end{bmatrix}
	\begin{bmatrix}
	   \Delta \Vec{\sigma^\res}\\ \Delta \Vec{\sigma^\well} \\ \Delta \Vec{p^\res} \\ \Delta \Vec{p^\well} \\ \Delta \Vec{s}
	\end{bmatrix}
	&=
	-
	\begin{bmatrix}
	   \Vec{r}^{\sigma^\res} \\ \Vec{r}^{\sigma^\well} \\ \Vec{r}^{p^\res} \\ \Vec{r}^{p^\well} \\ \Vec{r}^{s}
	\end{bmatrix},
	\label{eq:jac_system}
\end{align}
We solve \eqref{eq:jac_system} differently depending on which level is being solved.

\subsection{Reduction of fine-level Jacobian}
On the fine level, because the sub-matrix of the first $2\times2$ blocks of \eqref{eq:jac_system} is diagonal, \eqref{eq:jac_system} can easily be reduced to a system involving only $\Delta\Vec{p}^\res$, $\Delta\Vec{p}^\well$ and $\Delta\Vec{s}$ (the primal formulation) by eliminating $\Delta\Vec{\sigma}^\res$ and $\Delta\Vec{\sigma}^\well$.
\begin{align}
    \begin{bmatrix}
		0 & 0 & 0 \\
		0 & 0 & 0 \\
		0 & 0 & \Mat{A}_{ss} \\
	\end{bmatrix}
	-
	\begin{bmatrix}
		\Mat{A}_{p^\res\sigma^\res} & \Mat{A}_{p^r\sigma^\well} \\
		0 & \Mat{A}_{p^\well\sigma^\well} \\
		\Mat{A}_{s\sigma^\res} & \Mat{A}_{s\sigma^\well}\\
	\end{bmatrix}
	\begin{bmatrix}
		\Mat{A}_{\sigma^\res\sigma^\res} & \Mat{A}_{\sigma^\res\sigma^\well}  \\
		\Mat{A}_{\sigma^\well\sigma^\res} & \Mat{A}_{\sigma^\well\sigma^\well}  \\
	\end{bmatrix}^{-1}
	\begin{bmatrix}
		\Mat{A}_{\sigma^\res p^r} & 0 & \Mat{A}_{\sigma^\res s} \\
		\Mat{A}_{\sigma^\well p^\res} & \Mat{A}_{\sigma^\well p^\well} & \Mat{A}_{\sigma^\well s} \\
	\end{bmatrix}
	\label{eq:jac_system_primal}
\end{align}

\subsection{Reduction of coarse-level Jacobian}

On coarse levels, the sub-matrix of the first $2\times2$ blocks of \eqref{eq:jac_system} is not diagonal.
Therefore, instead of doing a direct inversion, we apply the idea of algebraic hybridization \cite{lee17, dobrev18} to eliminate the fluxes in \eqref{eq:jac_system}. 
More precisely, for each flux associated with an interior face in the computational grid, we introduce two one-sided fluxes, one for each of the grid cells sharing the face.
Both of the one-sided fluxes approximate the original flux on the face, and their continuity is enforced weakly by a Lagrange multiplier $\lambda_\varepsilon$ (also known as the face pressure in the context of subsurface flow).
We replace the two-point flux approximation \eqref{eq:TPFA_flux} by 
%
\begin{equation}
  \begin{split}
    \left( \frac{1}{\lambda(s_K) \overline{\Upsilon}_{K,\varepsilon}} \widehat{\sigma}_{\varepsilon, K} \right)  -  (p_{K} - \lambda_\varepsilon) & = 0,  \\
    \left(\frac{1}{\lambda(s_L) \overline{\Upsilon}_{L,\varepsilon}}   \widehat{\sigma}_{\varepsilon, L} \right)  -  (\lambda_\varepsilon - p_{L}) & = 0,  \\
    \widehat{\sigma}_{\varepsilon, K} - \widehat{\sigma}_{\varepsilon, L} & = 0.
  \end{split}
  \label{eq:TPFA_one_sided_flux}
\end{equation}
%
Note that the splitting of the weights for $\widehat{\sigma}_{\varepsilon, K}$ and $\widehat{\sigma}_{\varepsilon, L}$ naturally follows from the two-point flux approximation \eqref{eq:TPFA_flux}.
For flux equation \eqref{eq:total_volumetric_rate} associated with each well perforation, we want to use a similar idea so that our final system has a consistent structure.
However, unlike \eqref{eq:TPFA_flux}, the weight for $\sigma_{i}^\well$ is not the sum of the respective weights from two adjacent cells.
Therefore, we artificially create an algebraic splitting of the weight, and replace \eqref{eq:total_volumetric_rate} by:
\begin{equation}
  \begin{split}
    \frac{\alpha}{\lambda(s_{K(i)})WI_i} \widehat{\sigma}_{i, \well}^\well - \left( p^\well - \lambda^\well_i \right) & = 0 \\
    \frac{1-\alpha}{\lambda(s_{K(i)})WI_i} \widehat{\sigma}_{i, K(i)}^\well - \left( \lambda^\well_i - p^\res_{K(i)} \right) & = 0 \\
    \widehat{\sigma}_{i_\well}^\well - \widehat{\sigma}_{i, K(i)}^\well = 0,
  \end{split}
  \label{eq:total_volumetric_rate_hybrid}
\end{equation}
where $\alpha \in (0, 1)$ is a user constant, $\widehat{\sigma}_{i, \well}^\well$ and $\widehat{\sigma}_{i, K(i)}^\well$ are artificial fluxes replacing $\sigma_{i}^\well$, and $\lambda^\well_i$ is the associated Lagrange multiplier. 
In our numerical examples, we took $\alpha = 1/2$.
It is straight forward to see that \eqref{eq:TPFA_one_sided_flux} is equivalent to \eqref{eq:TPFA_flux}, and \eqref{eq:total_volumetric_rate_hybrid} is equivalent to \eqref{eq:total_volumetric_rate}
The discrete problem with one-sided fluxes is
\begin{align}
    \begin{bmatrix}
		\widehat{\Mat{A}}_{\sigma^\res\sigma^\res} & \widehat{\Mat{A}}_{\sigma^\res\sigma^\well} & \widehat{\Mat{A}}_{\sigma^\res p^\res} &   & \Mat{A}_{\sigma^\res\lambda^\res} &  & \widehat{\Mat{A}}_{\sigma^\res s} \\
		\widehat{\Mat{A}}_{\sigma^\well\sigma^\res} & \widehat{\Mat{A}}_{\sigma^\well\sigma^\well} & \widehat{\Mat{A}}_{\sigma^\well p^\res} & \widehat{\Mat{A}}_{\sigma^\well p^\well} &  & \Mat{A}_{\sigma^\well\lambda^\well} & \widehat{\Mat{A}}_{\sigma^\well s} \\
		\widehat{\Mat{A}}_{p^\res\sigma^\res} & \widehat{\Mat{A}}_{p^r\sigma^\well} &  &  &  & \\
		 & \widehat{\Mat{A}}_{p^\well\sigma^\well} &  &  &  & \\
		\Mat{A}_{\lambda^\res\sigma^\res} &  &  &  &  & \\
		 & \Mat{A}_{\lambda^\well\sigma^\well} &  &  &  & \\
		\widehat{\Mat{A}}_{s\sigma^\res} & \widehat{\Mat{A}}_{s\sigma^\well} &  &  &  &  & \Mat{A}_{ss}
	\end{bmatrix}
	\begin{bmatrix}
	   \Delta \widehat{\Vec{\sigma}}^\res\\ \Delta \widehat{\Vec{\sigma}}^\well \\ \Delta \Vec{p}^\res \\ \Delta \Vec{p}^\well \\ \Delta \Vec{\lambda}^\res \\ \Delta \Vec{\lambda}^\well \\ \Delta \Vec{s}
	\end{bmatrix}
	&=
	-
	\begin{bmatrix}
	   \Vec{r}^{\sigma^\res} \\ \Vec{r}^{\sigma^\well} \\ \Vec{r}^{p^\res} \\ \Vec{r}^{p^\well} \\ \Vec{r}^{\lambda^\res} \\ \Vec{r}^{\lambda^\well} \\ \Vec{r}^{s}
	\end{bmatrix},
	\label{eq:jac_system_hybrid}
\end{align}
In view of \eqref{eq:TPFA_one_sided_flux} and \eqref{eq:total_volumetric_rate_hybrid}, the first $4\times4$ blocks of \eqref{eq:jac_system_hybrid} can be rearranged into a block-diagonal matrix with invertible small blocks, each associated with either a grid cell or a well. Hence, \eqref{eq:jac_system_hybrid} can be efficiently reduced to a system for $\Delta\Vec{\lambda}^\res$, $\Delta\Vec{\lambda}^\well$ and $\Delta\Vec{s}$:
\begin{align}
    \begin{bmatrix}
		0 & 0 & 0 \\
		0 & 0 & 0 \\
		0 & 0 & \Mat{A}_{ss}
	\end{bmatrix}
	-
    \begin{bmatrix}
		\Mat{A}_{\lambda^\res\sigma^\res} & 0 & 0 & 0 \\
		0 & \Mat{A}_{\lambda^\well\sigma^\well} & 0 & 0 \\
		\widehat{\Mat{A}}_{s\sigma^\res} & \widehat{\Mat{A}}_{s\sigma^\well} & 0 & 0
	\end{bmatrix}
	\begin{bmatrix}
		\widehat{\Mat{A}}_{\sigma^\res\sigma^\res} & \widehat{\Mat{A}}_{\sigma^\res\sigma^\well} & \widehat{\Mat{A}}_{\sigma^\res p^\res} &   \\
		\widehat{\Mat{A}}_{\sigma^\well\sigma^\res} & \widehat{\Mat{A}}_{\sigma^\well\sigma^\well} & \widehat{\Mat{A}}_{\sigma^\well p^\res} & \widehat{\Mat{A}}_{\sigma^\well p^\well}  \\
		\widehat{\Mat{A}}_{p^\res\sigma^\res} & \widehat{\Mat{A}}_{p^r\sigma^\well} &  &  \\
		 & \widehat{\Mat{A}}_{p^\well\sigma^\well} &  & \\
	\end{bmatrix}^{-1}
    \begin{bmatrix}
		\Mat{A}_{\sigma^\res\lambda^\res} & 0 & \widehat{\Mat{A}}_{\sigma^\res s} \\
		0 & \Mat{A}_{\sigma^\well\lambda^\well} & \widehat{\Mat{A}}_{\sigma^\well s} \\
		0 & 0 & 0 \\
		0 & 0 & 0 \\
	\end{bmatrix}
	\label{eq:jac_system_hybrid_reduced}
\end{align}

\begin{rmk}
The physical meaning of the artificially created Lagrange multiplier is the pressure value at some convex combination of the well perforation and the cell center of the perforated cell. 
\end{rmk}

\begin{rmk}
Note that \eqref{eq:jac_system_primal} and \eqref{eq:jac_system_hybrid_reduced} have similar structures, so we solve them using the same linear solver.
For \eqref{eq:jac_system_primal}, we consider $\Delta\Vec{p}^\res$, $\Delta\Vec{p}^\well$ as one pressure block. 
Similarly, we group $\Delta\Vec{\lambda}^\res$, $\Delta\Vec{\lambda}^\well$ as one ``face pressure" block for \eqref{eq:jac_system_hybrid_reduced}. 
Then, we solve systems \eqref{eq:jac_system_primal} and \eqref{eq:jac_system_hybrid_reduced} by GMRES preconditioned by the CPR-type preconditioner in Section~3.4.3 of \cite{fas-two-phase}. 
\end{rmk}

\section{Numerical examples}\label{sec:numerical_examples}

In this section, we present three challenging numerical examples to demonstrate that the proposed methods can efficiently handle complex reservoir flow problems coupled with multi-perforation wells. 
In Section~\ref{sec:synthetic_numerical_example}, we use a synthetic three-layer test case with 25 wells to illustrate the applicability of FAS to cases with a large number of perforated cells.
Then, in Sections~\ref{sec:egg_model} and \ref{sec:saigup_model}, we include multi-perforation wells in the test cases first used in \cite{fas-two-phase} to assess the robustness of the well treatment in FAS on realistic test cases with sharp saturation fronts and high geometric complexity.
The simulation parameters are summarized in Table~\ref{tab:parameters} and the problem sizes are given in Table~\ref{tab:problem_sizes}.

In the three cases presented below, coarse cell aggregates are generated using METIS \cite{karypis1998fast} without using information on the intrinsic structure of the mesh. 
The coarsening procedure applied to the near-well regions is detailed in Section~\ref{sec:aggregation_around_wells}.
The coarsening factor $\beta$ is computed using the average aggregate size on each level.
We use a time stepping strategy in which the time step size is multiplied by a factor $\nu >1$ at every step:
\begin{equation}
\Delta t_{m} = \nu \Delta t_{m-1}, \quad m \geq 1. \label{time_stepping}
\end{equation}
This aggressive time stepping scheme is chosen to increase the problem difficulty and test the robustness of the nonlinear solvers.
For each time step, we report the largest Courant-Friedrichs-Lewy (CFL) number \cite{cao2002development} observed in the mesh.
We use a tolerance of $10^{-6}$ on the numerical residual to check the convergence of the nonlinear solvers.
The maximum number of nonlinear iterations is set to 10 on the coarsest level, and to 1 on all other levels.
We use a standard local saturation chopping strategy \cite{younis2011modern} to improve nonlinear convergence.
Specifically, we force the saturations to remain in [0,1] after each fine-level FAS update and after each single-level Newton update.
After the coarse-level FAS updates, we do not use local saturation chopping and instead we extend the mobility functions with constant values outside [0,1].

The discrete problems are generated using our own implementation based on MFEM \cite{mfem} of the finite-volume scheme described in Section \ref{sec:FV_tpfa}.
The multilevel spectral coarsening is performed with smoothG \cite{smoothg}, and the visualization is generated with ParaView. 

All the experiments were performed on a single node (Intel Xeon E5-2695 v4 @ 2.10 GHz) of the cluster PASCAL at Lawrence Livermore National Laboratory, which has 18 cores and 256 GB of memory.

\begin{table}[htbp]
	\caption{Parameter values used for the numerical examples. For the three test cases, the times are reported in total pore volume injected (PVI), which is the ratio of the injected wetting-phase volume over the total pore volume of the reservoir.}
	\label{tab:parameters}
        \small
        \centering
        \begin{tabular}{lllllllll}
          \toprule
          Symbol & Parameter & Units &
          \begin{tabular}{@{}c@{}}Synthetic \\ lognormal \end{tabular} &         
          Egg & 
          SAIGUP\\ 
          \midrule
          $s^0$        &  Initial wetting-phase saturation    & [-]           & 0 & 0 & 0\\
          $\mu_w$      &  Wetting-phase viscosity             & [Pa$\cdot$s]        & $10^{-3}$ & $10^{-3}$ & $10^{-3}$ \\
          $\mu_{nw}$   &  Non-wetting phase viscosity         & [Pa$\cdot$s]        & $5.0 \times 10^{-3}$ & $5.0 \times 10^{-3}$ & $5.0 \times 10^{-3}$ \\          
          $\gamma$     &  Relative permeability exponent      & [-]           & 2 & 2, 3, or 4 & 2 \\
          $\lambda_{\alpha}$&  Phase mobility ($\alpha \in \{w,nw\}$)         & [Pa$^{-1}$$\cdot$s$^{-1}$] &  $s^{\gamma}_{\alpha}/\mu_{\alpha}$ &   $s^{\gamma}_{\alpha}/\mu_{\alpha}$ & $s^{\gamma}_{\alpha}/\mu_{\alpha}$ \\
          \midrule 
          $n^I$        &  Number of injectors                 & [-]           & 13 & 8 & 5 \\
          $n^P$        &  Number of producers                 & [-]           & 12 & 4 & 5 \\
          $q_w^I$      &  Wetting-phase injection rate        & [m$^3$$\cdot$s$^{-1}$] & $3 \times 10^{-5}$ & $1.3 \times 10^{-4}$ & $5 \times 10^{-2}$ \\
          $p_{bh}$     &  Bottomhole pressure                 & [Pa]          & $10^6$ & $10^6$ & $10^6$ \\
          \midrule
          $\Delta t_0$ &  Initial time step                   & [PVI]         & $3.0 \times 10^{-5}$ & $9.3 \times 10^{-5}$ & $1.0 \times 10^{-4}$ \\
          $\nu$        &  Time step increase factor           & [-]           & 2 & 2 & 2 \\
          $T_{f}$      &  Final time                          & [PVI]         & $6.6 \times 10^{-2}$ & $4.7 \times 10^{-2}$ & $5.3 \times 10^{-2}$ \\
          \bottomrule
        \end{tabular}
\end{table}


%

\begin{table}[htbp]
	\caption{Problem sizes in the numerical examples.}
	\label{tab:problem_sizes}
        \small
        \centering
        \begin{tabular}{lrrrrrr}
          \toprule
          &
          \begin{tabular}{@{}c@{}}Synthetic \\ lognormal \end{tabular} &
          \begin{tabular}{@{}c@{}}Egg \\ (refined) 
          \end{tabular} &          
          \begin{tabular}{@{}c@{}}SAIGUP \\ (refined) \end{tabular}\\
          \midrule
          $|\mathcal{T}|$               & 306,030 & 148,424 & 629,760  \\
          $|\mathcal{E}|$               & 934,351 & 431,092 & 1,912,471 \\
          
          Number of unknowns            & 1240,381 & 579,516 & 2,542,231 \\
          \bottomrule          
        \end{tabular}
\end{table}

\subsection{Synthetic example}
\label{sec:synthetic_numerical_example}

\begin{figure}[h!]
  \vspace{5mm}
  \begin{subfigure}[b]{0.45\textwidth}
    \centering
    \includegraphics[scale=.18,clip,trim=0 80 0 0]{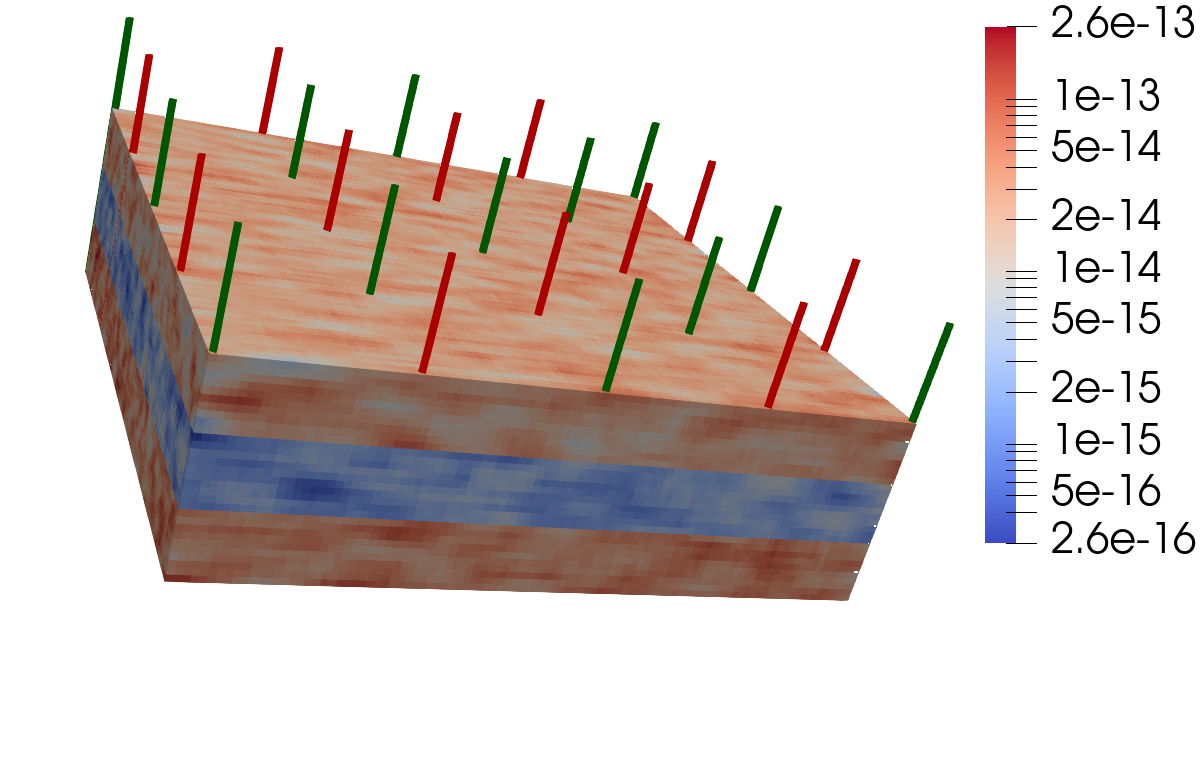}
    \caption{Permeability}
    \label{fig:permeability_synthetic}
  \end{subfigure}
  \hspace{8mm}
  \begin{subfigure}[b]{0.45\textwidth}
    \centering
    \includegraphics[scale=.18,clip,trim=0 80 0 0]{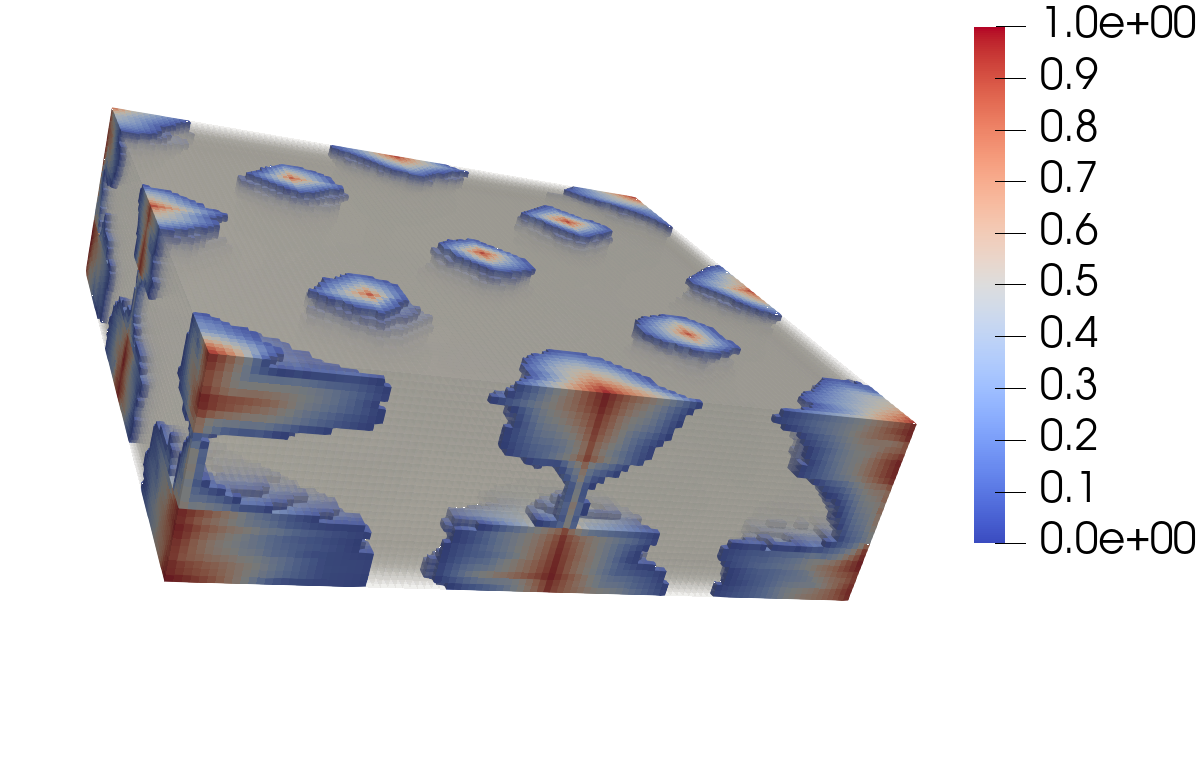}
    \caption{Saturation at final time}
    \label{fig:saturation_synthetic}
  \end{subfigure}
\caption{A 25-well synthetic example. Location of injectors (respectively producers) are indicated by green (respectively red) bars. For better visualization, the cells are rescaled to have the same length in all dimensions.}
\end{figure}

We first consider a synthetic example with a large number of wells and a large number of perforated cells to demonstrate the robustness and efficiency of the proposed FAS algorithm in this configuration.
We consider a three-dimensional domain of size 307.848 $\times$ 307.848 $\times$ 9.144 m$^3$ discretized with a $101 \times 101 \times 30$ Cartesian mesh consisting of 306,030 cells (each has size 3.048 $\times$ 3.048 $\times$ 0.3048 m$^3$).
The domain consists of three horizontal layers illustrated in Fig.~\ref{fig:permeability_synthetic}. 
On each layer, permeability follows a log-normal distribution provided by MRST \cite{lie19}.
Specifically, the permeability is generated by first calling the MRST function {\tt logNormalLayers([1010, 1010, 30], [7 1 5], 'indices', [1 11 21 31])}, and then it is rescaled to range between $2.6\times10^{-16}$ and $2.6\times10^{-13}$.
The permeability of the middle layer is approximately three orders of magnitude smaller than that of the top and bottom layers.
We place 25 wells with a spacing of 200 m in each direction as indicated in Fig.~\ref{fig:permeability_synthetic}.
Each well perforates a column of 30 cells from the top to the bottom of the reservoir ($n_{perf} = 30$), for a total of 750 cells perforated by a well in the reservoir.
The Peaceman indices $WI$ of equation \eqref{eq:phase_volumetric_rate} are computed using the standard formula given in \cite{Pea78}.
The fluid properties, well controls, and timestepping strategy are documented in Table~\ref{tab:parameters}.
The wetting-phase saturation map at the end of the simulation is shown in Fig.~\ref{fig:saturation_synthetic}.

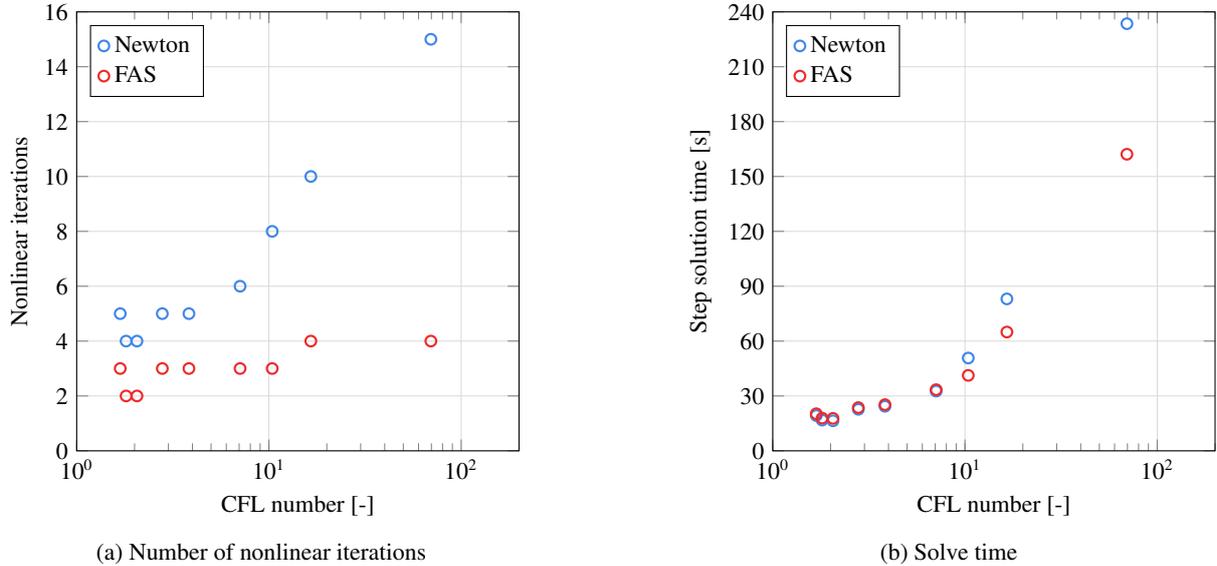
\begin{figure}[h!]
  \small
  \centering
  \begin{subfigure}[b]{0.45\textwidth}
    \centering
     \begin{tikzpicture}
    \begin{semilogxaxis}[
                  width=\textwidth,
                  height=\textwidth,
                  grid = major,
                  major grid style={very thin,draw=gray!30},
                  xmin=1,xmax=200,
                  ymin=0,ymax=16,
                  xlabel={CFL number [-]},
                  ylabel={Nonlinear iterations},
                  xtick={0.01, 0.1, 1, 10, 100},
                  ytick={0,2,...,20},
                  ylabel near ticks,
                  xlabel near ticks,
                  legend style={font=\small},
                  tick label style={font=\small},
                  label style={font=\small},
                  legend cell align={left},
                  legend pos=north west,
                  ]
      \addplot [mark=o, skyblue1, only marks, thick] table[x=CFL, y=nonlinear_iter, col sep=comma] {./data/lc/all_step_mrst_lognormal_ir3e-05_ro2_1lvl_cf128_4layers_scale1e6.csv};
      \addlegendentry{Newton};

      \addplot [mark=o, scarletred1, only marks, thick] table[x=CFL, y=nonlinear_iter, col sep=comma] {./data/lc/all_step_mrst_lognormal_ir3e-05_ro2_3lvl_cf128_4layers_scale1e6_merge_all_only1.csv};
      \addlegendentry{FAS};

            
    \end{semilogxaxis}
  \end{tikzpicture}
  
    \caption{Number of nonlinear iterations}
    \label{fig:nonlinear_convergence_synthetic}
  \end{subfigure}
  \hfill
  \begin{subfigure}[b]{0.45\textwidth}
    \centering
     \begin{tikzpicture}
    \begin{semilogxaxis}[
                  width=\textwidth,
                  height=\textwidth,
                  grid = major,
                  major grid style={very thin,draw=gray!30},
                  xmin=1,xmax=200,
                  ymin=0,ymax=240,
                  xlabel={CFL number [-]},
                  ylabel={Step solution time [s]},
                  xtick={0.01, 0.1, 1, 10, 100, 1000, 10000},
                  ytick={0,30,...,240},
                  ylabel near ticks,
                  xlabel near ticks,
                  legend style={font=\small},
                  tick label style={font=\small},
                  label style={font=\small},
                  legend cell align={left},
                  legend pos=north west,
                  ]
      \addplot [mark=o, skyblue1, only marks, thick] table[x=CFL, y=step_time, col sep=comma] {./data/lc/all_step_mrst_lognormal_ir3e-05_ro2_1lvl_cf128_4layers_scale1e6.csv};
      \addlegendentry{Newton};

      \addplot [mark=o, scarletred1, only marks, thick] table[x=CFL, y=step_time, col sep=comma] {./data/lc/all_step_mrst_lognormal_ir3e-05_ro2_3lvl_cf128_4layers_scale1e6_merge_all_only1.csv};
      \addlegendentry{FAS};

    \end{semilogxaxis}
  \end{tikzpicture}
  
    \caption{Solve time}
    \label{fig:step_solution_time_synthetic}
  \end{subfigure}

  \caption{Number of nonlinear iterations per time step as a function of CFL number [-] for the unfavorable end-point mobility ratio in the synthetic test case (MRST lognormal perm and 25 wells). FAS has 4 levels, with coarsening factor = 64.}
  \label{fig:mrst5_CFL_iterations_smear}
\end{figure}

Figure~\ref{fig:nonlinear_convergence_synthetic} compares the number of nonlinear iterations performed by FAS and single-level Newton during the simulation.
We make a distinction between two regimes based on the CFL number.
For the time steps corresponding to small CFL numbers (smaller than four), the nonlinear behavior of the two solution strategies is stable.
During this phase, FAS converges slightly faster than single-level Newton.
For the time steps corresponding to larger CFL numbers (larger than four), the number of nonlinear iterations required by single-level Newton rapidly deteriorates while that of FAS remains almost constant.
For the last time step of the simulation, single-level Newton takes 15 iterations to converge, compared to just four with FAS.
This example confirms that the algorithm proposed in this article  can achieve the excellent nonlinear behavior observed previously in \cite{fas-two-phase} but now in the presence of multi-perforation wells.

In terms of step solution time, Fig.~\ref{fig:step_solution_time_synthetic} shows that the mild reduction in nonlinear iterations obtained with FAS during the first time steps is not sufficient to reduce the solution time. 
However, as the CFL number increases and the convergence of single-level Newton starts deteriorating, the reduction in solution time obtained by FAS becomes more significant.
The largest reduction is obtained for the last time step and is equal to 31\%.

\subsection{Refined Egg model}
\label{sec:egg_model}

\begin{figure}[h!]
  \vspace{8mm}
  \begin{subfigure}[b]{0.45\textwidth}
    \centering
    \includegraphics[scale=.18,clip,trim=0 80 0 0]{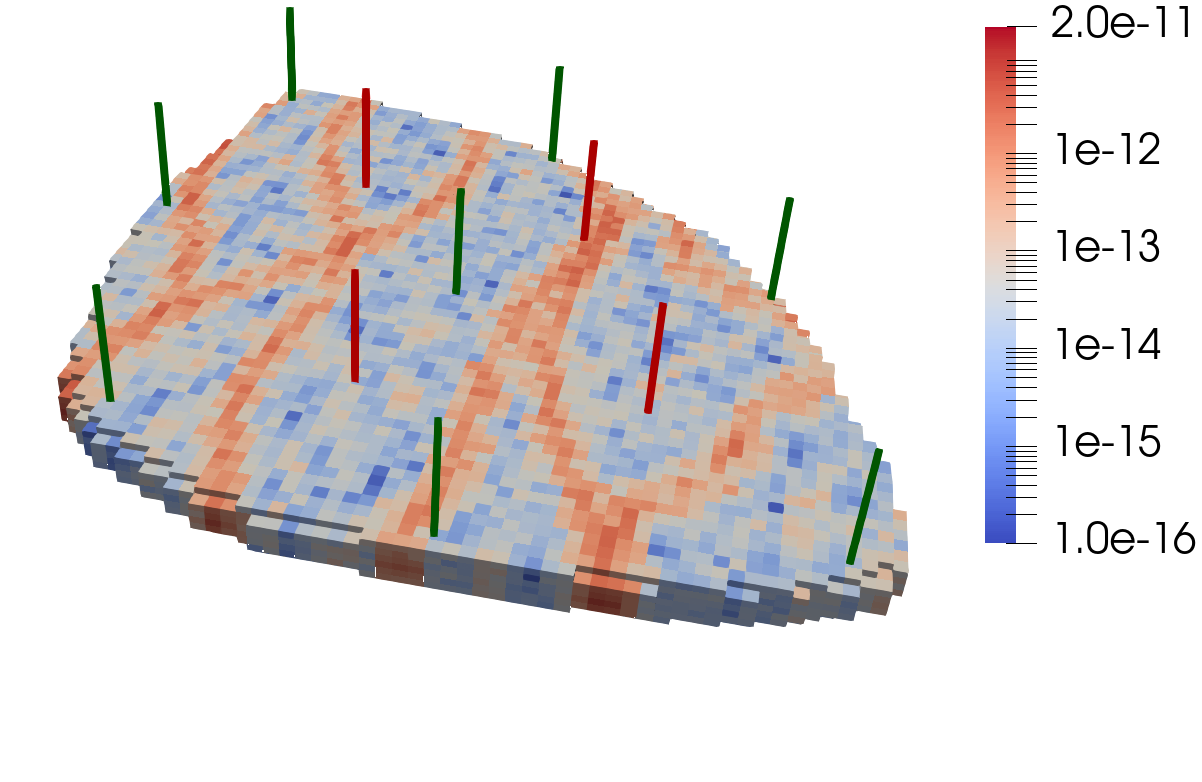}
    \caption{Permeability in x-direction.}
    \label{fig:perm-egg}
  \end{subfigure}
  \hspace{8mm}
  \begin{subfigure}[b]{0.45\textwidth}
    \centering
    \includegraphics[scale=.18,clip,trim=0 80 0 0]{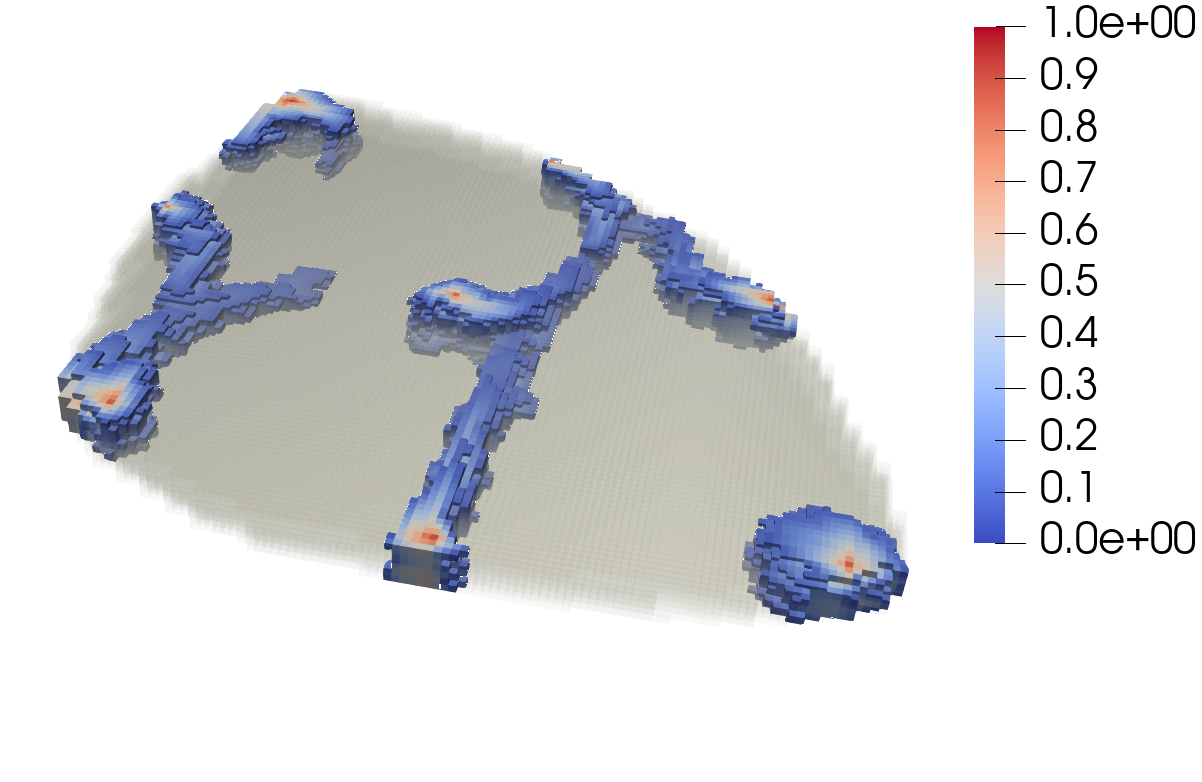}
    \caption{Saturation at final time}
    \label{fig:sat-egg}
  \end{subfigure}
\caption{A refined version of the Egg model. Location of injectors (respectively producers) are indicated by green (respectively red) bars.}
\end{figure}

In this section, we consider a refined version of the Egg model with a mesh consisting of 148,424 active cells, which is a $2 \times 2 \times 2$ regular refinement of the original mesh used in the published benchmark \cite{EggModel}.
We make the model highly heterogeneous by rescaling the permeability field and imposing a ratio of $2 \times 10^5$ between the largest and the smallest permeability in each direction.
As in the published benchmark, we use homogeneous porosity with $\Phi = 0.2$.
This example aims at demonstrating that FAS can preserve a robust nonlinear behavior in the presence of wells perforating the full model thickness with a channelized permeability field.
We follow the procedure employed in \cite{fas-two-phase} and assess the robustness of FAS when the strength of the nonlinearity increases.
This is done by increasing the relative permeability exponent from 2 to 4 in the mobility function definition.
The twelve wells (eight injectors and four producers) are placed using the $x-y$ coordinates provided in the published benchmark.
At these locations, each well perforates a column of 14 cells, and the number of perforation per well is 14 ($n_{perf} = 14$).
In total, 168 reservoir cells are perforated by a well in this test case.
The Peaceman indices $WI$ of equation \eqref{eq:phase_volumetric_rate} for each well perforation are computed using GEOS \cite{GEOS}.
The timestepping is the same as in the previous test case with $\nu = 2$.
%
Figure \ref{fig:sat-egg} shows the final wetting-phase saturation map.

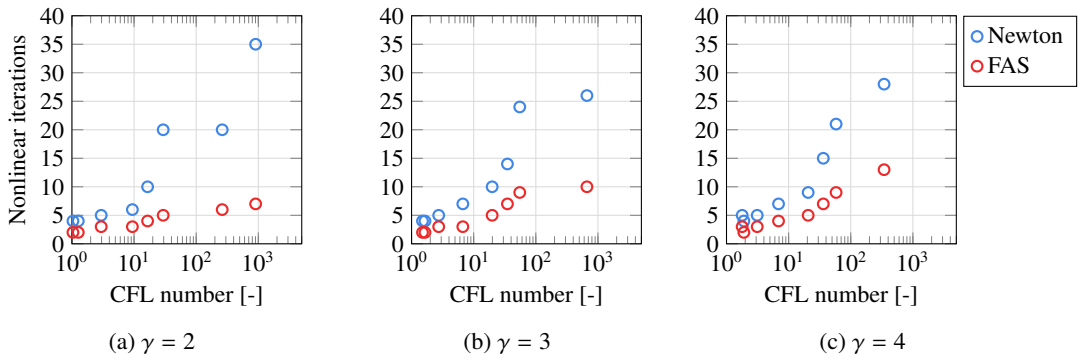
\begin{figure}[ht!]
  \footnotesize
  \begin{subfigure}[b]{0.28\textwidth}
    \centering
    \begin{tikzpicture}
      \begin{semilogxaxis}[
                  width=\textwidth,
                  height=\textwidth,
                  grid = major,
                  major grid style={very thin,draw=gray!30},
                  xmin=1,xmax=5000,
                  ymin=0,ymax=40,
                  xlabel={CFL number [-]},
                  ylabel={Nonlinear iterations},
                  xtick={0.01, 0.1, 1, 10, 100, 1000, 10000},
                  ytick={0,5,...,60},
                  ylabel near ticks,
                  xlabel near ticks,
                  legend style={font=\small},
                  tick label style={font=\small},
                  label style={font=\small},
                  legend cell align={left},
                  legend pos=north west,
                  ]
        \addplot [mark=o, skyblue1, only marks, thick] table[x=CFL, y=nonlinear_iter, col sep=comma] {./data/lc/all_step_egg_well_model_2x2x2_ir0.00013_ro2_1lvl_cf128_4layers_scale1e6.csv};
        \addplot [mark=o, scarletred1, only marks, thick] table[x=CFL, y=nonlinear_iter, col sep=comma] {./data/lc/all_step_egg_well_model_2x2x2_ir0.00013_ro2_3lvl_cf128_4layers_scale1e6_merge_all_only1.csv};
        
      \end{semilogxaxis}
    \end{tikzpicture}
    \caption{$\gamma = 2$}
  \end{subfigure}
  \begin{subfigure}[b]{0.28\textwidth}
    \centering
    \begin{tikzpicture}
      \begin{semilogxaxis}[
                  width=\textwidth,
                  height=\textwidth,
                  grid = major,
                  major grid style={very thin,draw=gray!30},
                  xmin=1,xmax=5000,
                  ymin=0,ymax=40,
                  xlabel={CFL number [-]},
                  xtick={0.01, 0.1, 1, 10, 100, 1000, 10000},
                  ytick={0,5,...,60},
                  ylabel near ticks,
                  xlabel near ticks,
                  legend style={font=\small},
                  tick label style={font=\small},
                  label style={font=\small},
                  legend cell align={left},
                  legend pos=north west,
                  ]
        \addplot [mark=o, skyblue1, only marks, thick] table[x=CFL, y=nonlinear_iter, col sep=comma] {./data/lc/all_step_egg_well_model_2x2x2_ir0.00013_ro3_1lvl_cf128_4layers_scale1e6.csv};
        \addplot [mark=o, scarletred1, only marks, thick] table[x=CFL, y=nonlinear_iter, col sep=comma] {./data/lc/all_step_egg_well_model_2x2x2_ir0.00013_ro3_3lvl_cf128_4layers_scale1e6_merge_all_only1.csv};
        
      \end{semilogxaxis}
    \end{tikzpicture}
    \caption{$\gamma = 3$}
  \end{subfigure}
  \begin{subfigure}[b]{0.28\textwidth}
    \centering
    \begin{tikzpicture}
      \begin{semilogxaxis}[
                  width=\textwidth,
                  height=\textwidth,
                  grid = major,
                  major grid style={very thin,draw=gray!30},
                  xmin=1,xmax=5000,
                  ymin=0,ymax=40,
                  xlabel={CFL number [-]},
                  xtick={0.01, 0.1, 1, 10, 100, 1000, 10000},
                  ytick={0,5,...,60},
                  ylabel near ticks,
                  xlabel near ticks,
                  legend style={font=\small},
                  tick label style={font=\small},
                  label style={font=\small},
                  legend cell align={left},
                  legend pos=outer north east,
                  ]
        \addplot [mark=o, skyblue1, only marks, thick] table[x=CFL, y=nonlinear_iter, col sep=comma] {./data/lc/all_step_egg_well_model_2x2x2_ir0.00013_ro4_1lvl_cf128_4layers_scale1e6.csv};
        \addlegendentry{Newton};
        \addplot [mark=o, scarletred1, only marks, thick] table[x=CFL, y=nonlinear_iter, col sep=comma] {./data/lc/all_step_egg_well_model_2x2x2_ir0.00013_ro4_3lvl_cf128_4layers_scale1e6_merge_all_only1.csv};
        \addlegendentry{FAS};
        
      \end{semilogxaxis}
    \end{tikzpicture}
    \caption{$\gamma = 4$}
  \end{subfigure}
\caption{Number of nonlinear iterations for the refined Egg model.}
   \label{fig:nonlinear_convergence_egg}
\end{figure}

We use a 3-level FAS with a coarsening factor $\beta = 32$. 
Figure~\ref{fig:well-aggregate} shows the shape of a cell aggregate in contact with a well to illustrate the aggregation strategy around wells described in Section~\ref{sec:aggregation_around_wells}.
The nonlinear behavior of FAS and single-level Newton with damping is presented in Fig.~\ref{fig:nonlinear_convergence_egg} for the three relative permeability exponents ($\gamma = 2, 3, 4$).
The three cases confirm the robustness of FAS with respect to time step size and nonlinearity of the relative permeability function.
At the beginning of the simulation, for CFL numbers smaller than 10, FAS achieves a mild reduction in the number of nonlinear iterations (between 2 and 4 nonlinear iterations saved per time step) compared to single-level Newton.
As the time step increases, the FAS nonlinear convergence remains robust, while that of single-level Newton deteriorates significantly. 
For the last time step of the simulation, which has a CFL number above 100 in the three cases, FAS only requires 7 ($\gamma = 2$), 10 ($\gamma = 3$), and 13 ($\gamma = 4$) nonlinear iterations, compared to 35 ($\gamma = 2$), 26 ($\gamma = 3$), and 28 ($\gamma = 4$) with single-level Newton.

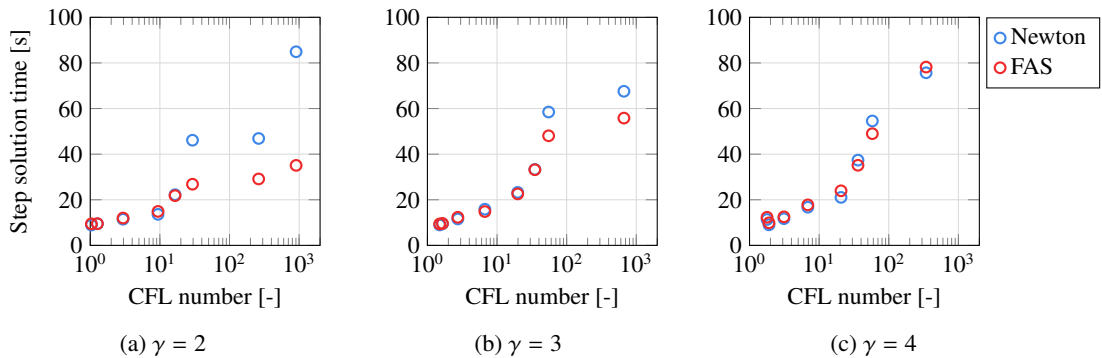
\begin{figure}[ht!]
  \footnotesize
  \begin{subfigure}[b]{0.28\textwidth}
    \centering
    \begin{tikzpicture}
      \begin{semilogxaxis}[
                  width=\textwidth,
                  height=\textwidth,
                  grid = major,
                  major grid style={very thin,draw=gray!30},
                  xmin=1,xmax=2000,
                  ymin=0.01,ymax=100,
                  xlabel={CFL number [-]},
                  ylabel={Step solution time [s]},
                  xtick={0.01, 0.1, 1, 10, 100, 1000, 10000},
                  ytick={0,20,...,120},
                  ylabel near ticks,
                  xlabel near ticks,
                  legend style={font=\small},
                  tick label style={font=\small},
                  label style={font=\small},
                  legend cell align={left},
                  legend pos=north west,
                  ]
        \addplot [mark=o, skyblue1, only marks, thick] table[x=CFL, y=step_time, col sep=comma] {./data/lc/all_step_egg_well_model_2x2x2_ir0.00013_ro2_1lvl_cf128_4layers_scale1e6.csv};
        \addplot [mark=o, scarletred1, only marks, thick] table[x=CFL, y=step_time, col sep=comma] {./data/lc/all_step_egg_well_model_2x2x2_ir0.00013_ro2_3lvl_cf128_4layers_scale1e6_merge_all_only1.csv};
            
      \end{semilogxaxis}
    \end{tikzpicture}
    \caption{$\gamma = 2$}
  \end{subfigure}
  \begin{subfigure}[b]{0.28\textwidth}
    \centering
    \begin{tikzpicture}
      \begin{semilogxaxis}[
                  width=\textwidth,
                  height=\textwidth,
                  grid = major,
                  major grid style={very thin,draw=gray!30},
                  xmin=1,xmax=2000,
                  ymin=0.01,ymax=100,
                  xlabel={CFL number [-]},
                  xtick={0.01, 0.1, 1, 10, 100, 1000, 10000},
                  ytick={0,20,...,120},
                  ylabel near ticks,
                  xlabel near ticks,
                  legend style={font=\small},
                  tick label style={font=\small},
                  label style={font=\small},
                  legend cell align={left},
                  legend pos=north west,
                  ]
        \addplot [mark=o, skyblue1, only marks, thick] table[x=CFL, y=step_time, col sep=comma] {./data/lc/all_step_egg_well_model_2x2x2_ir0.00013_ro3_1lvl_cf128_4layers_scale1e6.csv};
        \addplot [mark=o, scarletred1, only marks, thick] table[x=CFL, y=step_time, col sep=comma] {./data/lc/all_step_egg_well_model_2x2x2_ir0.00013_ro3_3lvl_cf128_4layers_scale1e6_merge_all_only1.csv};
            
      \end{semilogxaxis}
    \end{tikzpicture}
    \caption{$\gamma = 3$}
  \end{subfigure}
  \begin{subfigure}[b]{0.28\textwidth}
    \centering
    \begin{tikzpicture}
      \begin{semilogxaxis}[
                  width=\textwidth,
                  height=\textwidth,
                  grid = major,
                  major grid style={very thin,draw=gray!30},
                  xmin=1,xmax=2000,
                  ymin=0.01,ymax=100,
                  xlabel={CFL number [-]},
                  xtick={0.01, 0.1, 1, 10, 100, 1000, 10000},
                  ytick={0,20,...,120},
                  ylabel near ticks,
                  xlabel near ticks,
                  legend style={font=\small},
                  tick label style={font=\small},
                  label style={font=\small},
                  legend cell align={left},
                  legend pos=outer north east,
                  ]
        \addplot [mark=o, skyblue1, only marks, thick] table[x=CFL, y=step_time, col sep=comma] {./data/lc/all_step_egg_well_model_2x2x2_ir0.00013_ro4_1lvl_cf128_4layers_scale1e6.csv};
        \addlegendentry{Newton};
        \addplot [mark=o, scarletred1, only marks, thick] table[x=CFL, y=step_time, col sep=comma] {./data/lc/all_step_egg_well_model_2x2x2_ir0.00013_ro4_3lvl_cf128_4layers_scale1e6_merge_all_only1.csv};
        \addlegendentry{FAS};
            
      \end{semilogxaxis}
    \end{tikzpicture}
    \caption{$\gamma = 4$}
  \end{subfigure}
\caption{Solve time for the refined Egg model.}
\label{fig:step_solution_time_egg}
\end{figure}

Figure \ref{fig:step_solution_time_egg} summarizes the impact of the nonlinear behavior of the two solution strategies on the step solution time.
The mild reduction in the number of nonlinear iterations obtained with FAS for small CFL numbers is not sufficient to reduce the step solution time. 
However, for large CFL numbers, the more significant reduction in nonlinear iterations allows FAS to achieve a large reduction in solution time for $\gamma = 2$ and $\gamma = 3$ compared to single-level Newton.
Specifically, for the last time step, FAS reduces the solution time by 58\% and 17\% for $\gamma = 2$ and $\gamma = 3$, respectively.
For $\gamma = 4$, the solution time of FAS and single-level Newton remains comparable until the end of the simulation due to the large number of iterations required by the linear solver in the coarse FAS solves.

\subsection{Refined SAIGUP}
\label{sec:saigup_model}

\begin{figure}[h!]
  \centering
  \begin{subfigure}[b]{0.45\textwidth}
    \centering
    \includegraphics[scale=.19,clip,trim=0 80 0 0]{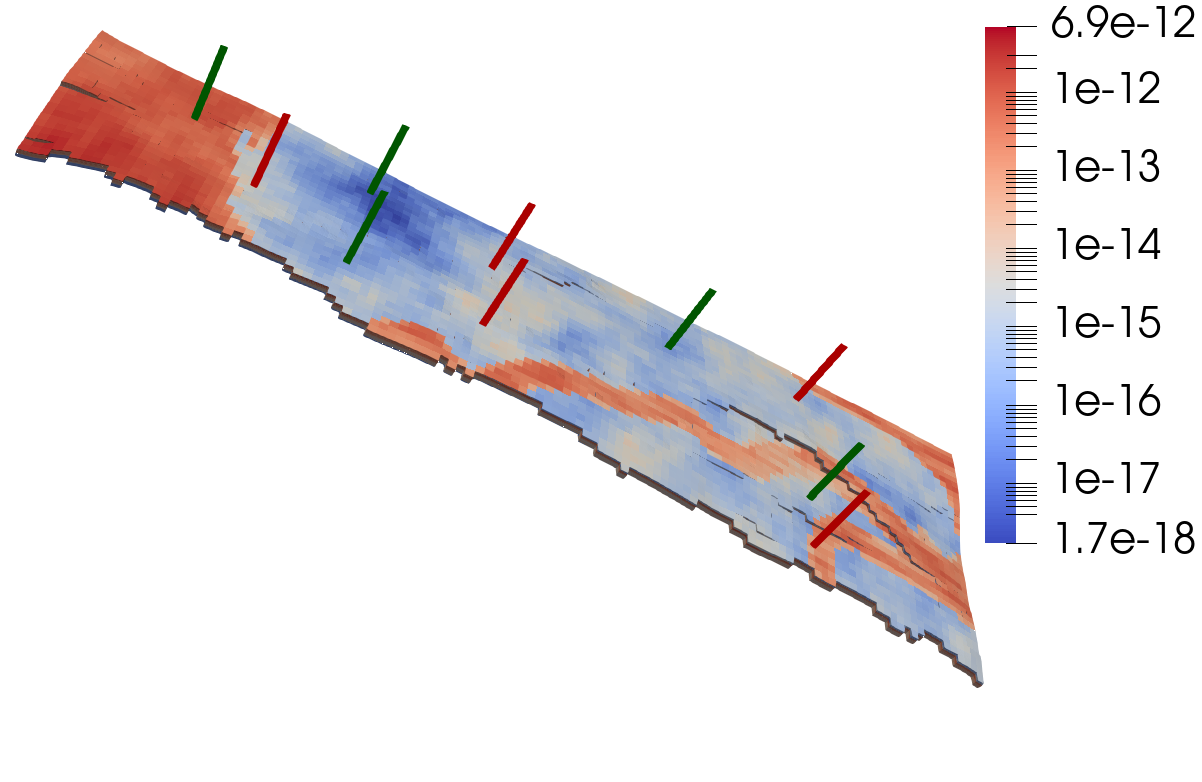}
    \caption{Permeability in x-direction}
    \label{fig:perm-x-saigup}
  \end{subfigure}
  \hspace{12mm}
  \begin{subfigure}[b]{0.45\textwidth}
    \centering
    \includegraphics[scale=.19,clip,trim=0 80 0 0]{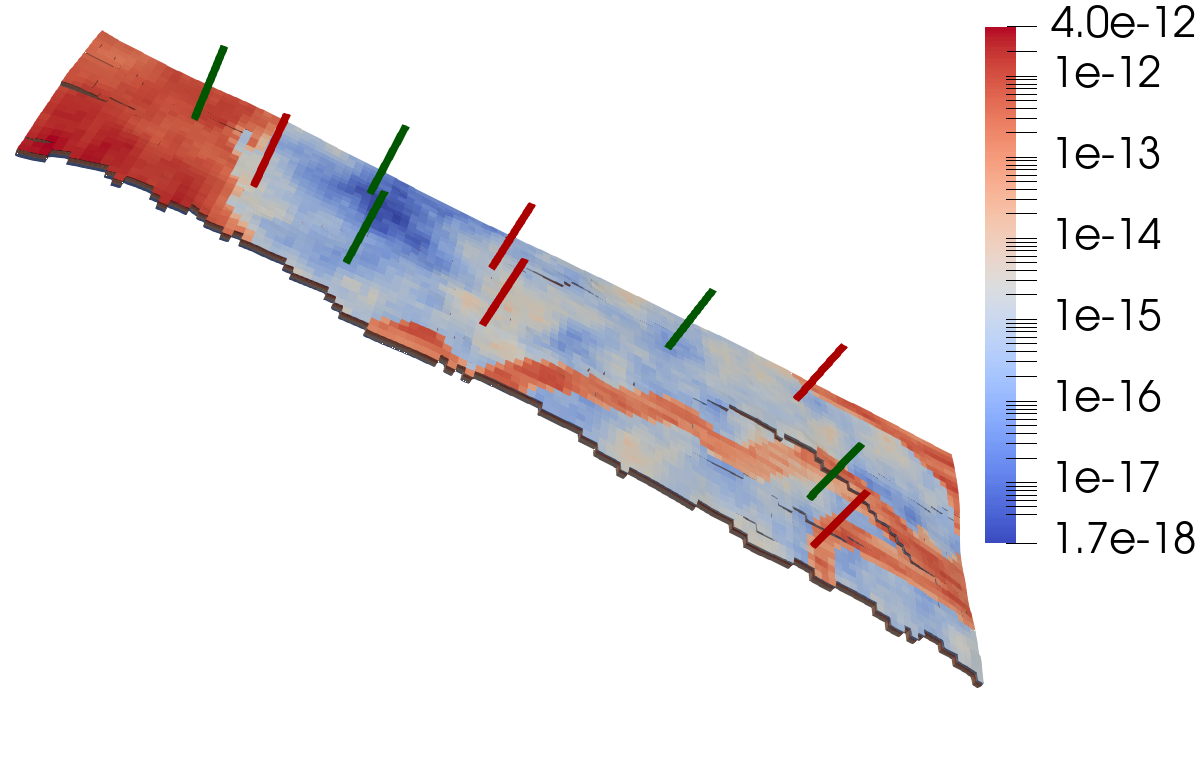}
    \caption{Permeability in y-direction}
    \label{fig:perm-y-saigup}
  \end{subfigure}
  \\
  \vspace{10mm}
  \begin{subfigure}[b]{0.45\textwidth}
    \centering
    \includegraphics[scale=.19,clip,trim=0 80 0 0]{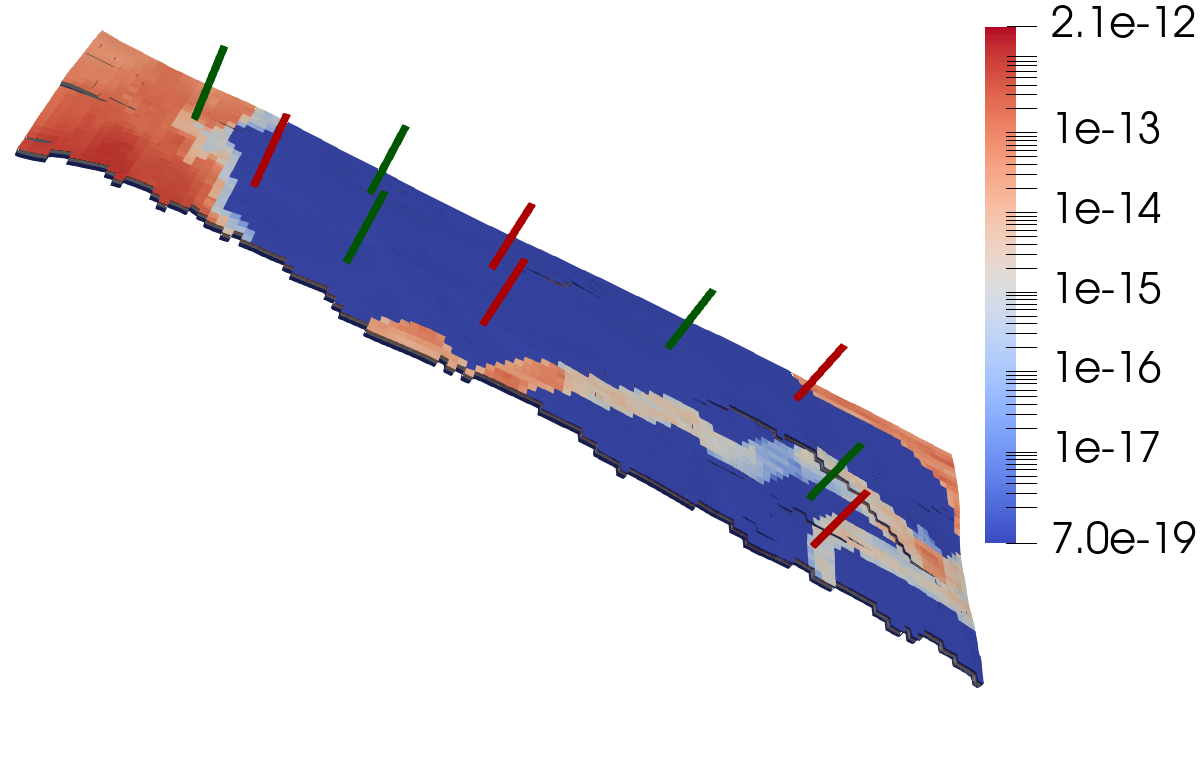}
    \caption{Permeability in z-direction}
    \label{fig:perm-z-saigup}
  \end{subfigure}
  \hspace{12mm}
  \begin{subfigure}[b]{0.45\textwidth}
    \centering
    \includegraphics[scale=.19,clip,trim=0 80 0 0]{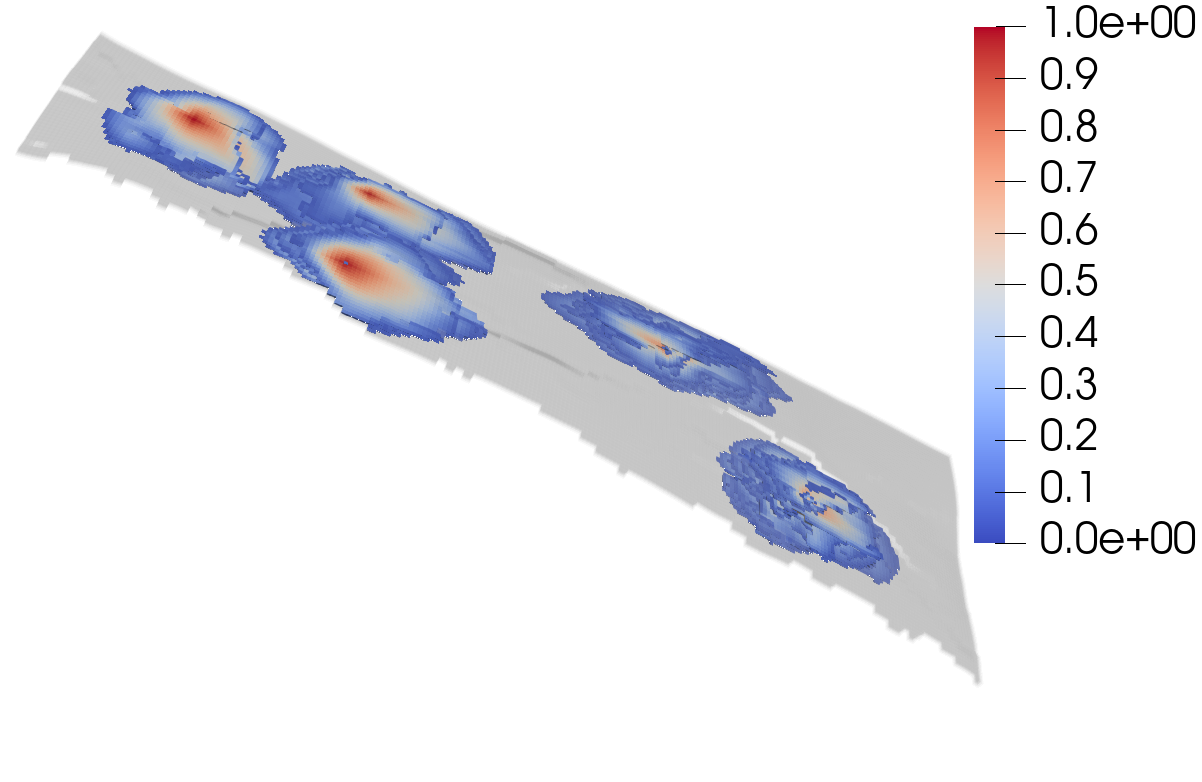}
    \caption{Saturation at final time}
    \label{fig:sat-saigup}
  \end{subfigure}
\caption{Permeability and saturation at final time of the refined version of the SAIGUP model. Location of injectors (respectively producers) are indicated by green (respectively red) bars.}
\end{figure}

This test case is derived from the SAIGUP model \cite{manzocchi2008sensitivity}. 
It is based on a regularly refined version of the original mesh and consists of 629,760 cells.
We use the permeability and porosity fields of the original model.
This example aims at demonstrating the robustness of FAS on a complex corner-point mesh with ten wells (five injectors and five producers) perforating the full domain thickness.
The wells are placed in the locations specified in the original problem.
Each well perforates 20 cells and the corresponding Peaceman indices are computed using MRST \cite{lie19}.
In total, 200 cells are perforated by a well.
We refer the reader to Table~\ref{tab:parameters} for additional details on the problem setup.
The final wetting-phase saturation maps are shown in Fig.~\ref{fig:sat-saigup}.

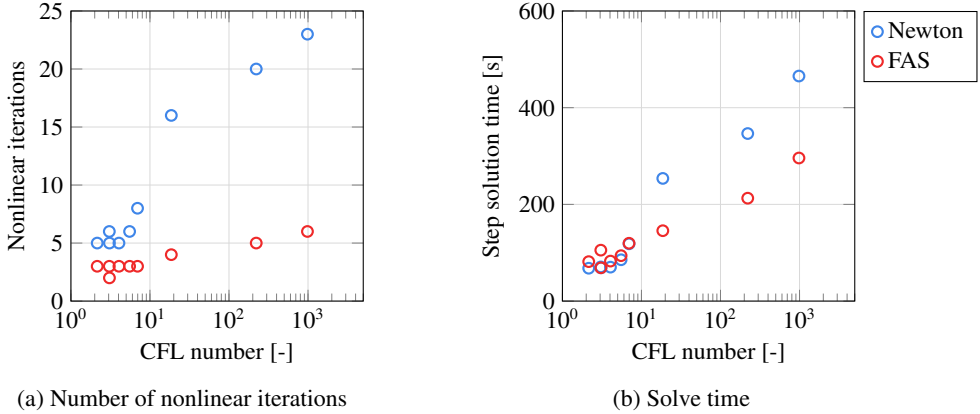
\begin{figure}[ht!]
  \footnotesize
  \begin{subfigure}[b]{0.33\textwidth}
    \centering
    \begin{tikzpicture}
      \begin{semilogxaxis}[
                  width=\textwidth,
                  height=\textwidth,
                  grid = major,
                  major grid style={very thin,draw=gray!30},
                  xmin=1,xmax=5000,
                  ymin=0,ymax=25,
                  xlabel={CFL number [-]},
                  ylabel={Nonlinear iterations},
                  xtick={0.01, 0.1, 1, 10, 100, 1000, 10000},
                  ytick={0,5,...,60},
                  ylabel near ticks,
                  xlabel near ticks,
                  legend style={font=\small},
                  tick label style={font=\small},
                  label style={font=\small},
                  legend cell align={left},
                  legend pos=outer north east,
                  ]
        \addplot [mark=o, skyblue1, only marks, thick] table[x=CFL, y=nonlinear_iter, col sep=comma] {./data/lc/all_step_saigup_well_model_2x2x2_ir0.05_ro2_1lvl_cf128_4layers_scale1e6.csv};

        \addplot [mark=o, scarletred1, only marks, thick] table[x=CFL, y=nonlinear_iter, col sep=comma] {./data/lc/all_step_saigup_well_model_2x2x2_ir0.05_ro2_4lvl_cf32_4layers_scale1e6_merge_all_only1.csv};

      \end{semilogxaxis}
    \end{tikzpicture}
    \caption{Number of nonlinear iterations}
    \label{fig:nonlinear_convergence_saigup}
  \end{subfigure}
  \hspace{10mm}
  \begin{subfigure}[b]{0.33\textwidth}
    \centering
    \begin{tikzpicture}
      \begin{semilogxaxis}[
                  width=\textwidth,
                  height=\textwidth,
                  grid = major,
                  major grid style={very thin,draw=gray!30},
                  xmin=1,xmax=5000,
                  ymin=0.01,ymax=600,
                  xlabel={CFL number [-]},
                  ylabel={Step solution time [s]},
                  xtick={0.01, 0.1, 1, 10, 100, 1000, 10000},
                  ytick={0,200,...,1500},
                  ylabel near ticks,
                  xlabel near ticks,
                  legend style={font=\small},
                  tick label style={font=\small},
                  label style={font=\small},
                  legend cell align={left},
                  legend pos=outer north east,
                  ]
        \addplot [mark=o, skyblue1, only marks, thick] table[x=CFL, y=step_time, col sep=comma] {./data/lc/all_step_saigup_well_model_2x2x2_ir0.05_ro2_1lvl_cf128_4layers_scale1e6.csv};
        \addlegendentry{Newton};

        \addplot [mark=o, scarletred1, only marks, thick] table[x=CFL, y=step_time, col sep=comma] {./data/lc/all_step_saigup_well_model_2x2x2_ir0.05_ro2_4lvl_cf32_4layers_scale1e6_merge_all_only1.csv};
        \addlegendentry{FAS};
        
            
      \end{semilogxaxis}
    \end{tikzpicture}
    \caption{Solve time}
    \label{fig:step_solution_time_saigup}
  \end{subfigure}
\caption{Solve time for the refined SAIGUP model.}
\end{figure}

The nonlinear behavior of four-level FAS and single-level Newton is documented in Fig.~\ref{fig:nonlinear_convergence_saigup}.
The coarsening factor from level 0 (fine level) to level 1 is 32, while that for the subsequent coarsening is 8.
The results are in agreement with the observations made on the Egg model in Section~\ref{sec:egg_model}.
The nonlinear behavior of FAS remains very stable throughout the simulation, including at the end of the simulation when CFL numbers reach large values.
In fact, for this test case, FAS can always converge in fewer than six nonlinear iterations.
The nonlinear behavior of single-level Newton deteriorates quickly as the CFL number grows, and for the last time step, 23 Newton iterations are needed to achieve convergence.
The step solution time of the two schemes is shown in Fig.~\ref{fig:step_solution_time_saigup}.
For the time steps corresponding to CFL numbers smaller than 10, the convergence acceleration obtained with FAS is too limited to result in a reduction in solution time compared to single-level Newton. 
For the last three time steps of the simulation, FAS reduces the step solution time by 43\%, 39\%, and 36\%.

To conclude this section, we evaluate the impact of the number of FAS levels on the total solution time in Table~\ref{tab:impact_of_fas_levels_saigup}.
This is done by running the refined SAIGUP problem described above with two, three, and four levels.

\begin{table}[ht]
\centering
\begin{tabular}{ ccccccc } 
\hline
Number of levels && 1 & 2 & 3 & 4  \\
\hline
Total solving time && 1553s & 1422s & 1222s & 1236s \\ 
\hline
\end{tabular}
\caption{Comparison of total solve time (in seconds) when using FAS with different number of levels to solve the refined SAIGUP problem. Note that the 1-level FAS is the standard Newton's method.}
\label{tab:impact_of_fas_levels_saigup}
\end{table}

\section{Concluding remarks}

In the current paper, we extend our previous work on the full approximation scheme nonlinear multigrid solver for two-phase flow and transport problems developed in \cite{fas-two-phase}, with a focus on the treatment of wells.
In particular, we consider a standard well model with multiple perforations.
The multigrid hierarchy is constructed using a multilevel coarsening algorithm for graph Laplacians, where the graph Laplacian is extracted from the discretization of the elliptic component of the model problem.
The coarsening algorithm is aggregation-based and algebraic, which allows the coarse grid to be flexible and adjustable.
To take advantage of this flexibility, a partitioning strategy exploiting the location of well cells is proposed to produce aggregations such that well cells are not in the vicinity of the aggregate interface.
As noted in the literature, such coarse grids are favourable because coarse basis functions have a much better representation of the radial flow near wells.
When solving the Jacobian system on coarse levels, we use a hybridization technique to transform the original discrete problem to one that involves only wetting-phase saturations and Lagrange multipliers.
For local problems involving a well perforation, an algebraic splitting is proposed so that the structure of all the local problems is consistent.
The final discrete system is aware of the unknowns associated with the wells, and those associated with the reservoir.
This piece of information can be further used to construct more efficient linear solvers by treating well unknowns separately, which is a subject of our ongoing research.
Also, while the discussion in this paper is based on a TPFA finite volume discretization, the proposed solver can be constructed in a similar manner if the underlying discretization is the mimetic finite difference method or the mixed finite element method.

In our numerical experiments, the proposed FAS solver is compared against the standard Newton method on one synthetic and two modified benchmark reservoir models.
We observe that the FAS solver in general outperforms Newton's method when the CFL number is larger than 10, while the two solvers have similar performance when the CFL number is smaller than 10.
For the refined SAIGUP model, which is unstructured and more realistic, FAS is 36\%-43\% faster than Newton's method in terms of solving time for the time steps with a large CFL number.
%
These nonlinear convergence improvements for large time steps make it easier for reservoir simulation practitioners to focus solely on accuracy considerations when selecting a time step size, without being limited by the nonlinear solver capabilities.
%

\section*{Acknowledgements}
\label{sec::acknow}
Funding was provided by TotalEnergies through the FC-MAELSTROM project.
Portions of this work were performed under the auspices of the U.S. Department of
Energy by Lawrence Livermore National Laboratory under Contract DE-AC52-07-NA27344 (LLNL-JRNL-851214).

\bibliography{authorList,journalAbbreviations,proceedingCollectionNames,publisherNames,references_NEW}

\end{document}

%% file: pics/well.tex
\begin{tikzpicture}
    \tikzstyle gridColor=[black!20]
    
    
    
    

    \draw[gridColor] (-4, -1.4) -- (2, -1.4) ;
    \draw[gridColor] (-4, -0.4) -- (2, -0.4) ;
    \draw[gridColor] (-4, 0.6) -- (2, 0.6) ;
    \draw[gridColor] (-4, 1.6) -- (2, 1.6) ;
    \draw[gridColor] (-3.5, 1.6) -- (-3.5, -1.4) ;
    
    \draw[gridColor, dashed] (-1.5, 1.1) -- (4.5, 1.1) ;
    \draw[gridColor, dashed] (-1.5, 2.1) -- (4.5, 2.1) ;
    \draw[gridColor, dashed] (-1.5, 3.1) -- (4.5, 3.1) ;
    \draw[gridColor] (-1.5, 4.1) -- (4.5, 4.1) ;
    \draw[gridColor, dashed] (-1, 4.1) -- (-1, 1.1) ;
    
    \draw[gridColor, dashed] (-3.5, -1.4) -- (-1, 1.1);
    \draw[gridColor, dashed] (-3.5, -0.4) -- (-1, 2.1) ;
    \draw[gridColor, dashed] (-3.5, 0.6) -- (-1, 3.1) ;
    \draw[gridColor] (-3.5, 1.6) -- (-1, 4.1) ;
    
    \draw[gridColor] (2, -1.4) -- (4.5, 1.1) ;
    \draw[gridColor] (2, -0.4) -- (4.5, 2.1) ;
    \draw[gridColor] (2, 0.6) -- (4.5, 3.1) ;
    \draw[gridColor] (2, 1.6) -- (4.5, 4.1) ;
    \draw[gridColor] (2, 1.6) -- (2, -1.4) ;
    \draw[gridColor] (4.5, 4.1) -- (4.5, 1.1) ;

    \tikzstyle wellColor=[plum1]
    \draw (0.5, 3)[wellColor, thick] ellipse (0.5 and 0.2);
    \draw (0.0, 2)[wellColor, thick] arc[start angle=180,end angle=360,x radius=0.5,y radius=0.2]; 
    \draw (0.0, 1)[wellColor, thick] arc[start angle=180,end angle=360,x radius=0.5,y radius=0.2]; 
    \draw (0.0, 0)[wellColor, thick] arc[start angle=180,end angle=360,x radius=0.5,y radius=0.2]; 
    \draw[wellColor, thick] (0, 0) -- (0, 3) [] {};
    \draw[wellColor, thick] (1, 0) -- (1, 3) [] {};
    
    \node[] at (2.5, 0.4) [] {\footnotesize Perforation 1};
    \node[] at (2.5, 1.4) [] {\footnotesize Perforation 2};
    \node[] at (2.5, 2.4) [] {\footnotesize Perforation 3};
    \draw[-stealth, opacity=0.5] (1.5, 0.4) -- (0.3, 0.4);
    \draw[-stealth, opacity=0.5] (1.5, 1.4) -- (0.3, 1.4);
    \draw[-stealth, opacity=0.5] (1.5, 2.4) -- (0.3, 2.4);
    
    \tikzstyle perfColor=[chocolate1]
    \draw[perfColor, thick] (-0.4, 0.35) -- (0.2, 0.35) [] {};
    \draw[perfColor, thick] (-0.4, 0.4) -- (0.2, 0.4) [] {};
    \draw[perfColor, thick] (-0.4, 1.35) -- (0.2, 1.35) [] {};
    \draw[perfColor, thick] (-0.4, 1.4) -- (0.2, 1.4) [] {};
    \draw[perfColor, thick] (-0.4, 2.35) -- (0.2, 2.35) [] {};
    \draw[perfColor, thick] (-0.4, 2.4) -- (0.2, 2.4) [] {};
    
    \draw [decorate, decoration = {calligraphic brace}, thick] (-0.6, 0.37) --  (-0.6, 1.33);
    \draw [decorate, decoration = {calligraphic brace}, thick] (-0.6, 1.4) --  (-0.6, 2.36);
    
    \node[] at (-1.3, 0.9) [] {\footnotesize$\triangle h_{1,2}$};
    \node[] at (-1.3, 1.9) [] {\footnotesize$\triangle h_{2,3}$};
\end{tikzpicture}

%% file: pics/well_dofs.tex
\begin{tikzpicture}
    \tikzstyle gridColor=[black!20]
    
    
    
    
    
    
    \draw[gridColor] (-4, -1.4) -- (2, -1.4) ;
    \draw[gridColor] (-4, -0.4) -- (2, -0.4) ;
    \draw[gridColor] (-4, 0.6) -- (2, 0.6) ;
    \draw[gridColor] (-4, 1.6) -- (2, 1.6) ;
    \draw[gridColor] (-3.5, 1.6) -- (-3.5, -1.4) ;
    
    \draw[gridColor, dashed] (-1.5, 1.1) -- (4.5, 1.1) ;
    \draw[gridColor, dashed] (-1.5, 2.1) -- (4.5, 2.1) ;
    \draw[gridColor, dashed] (-1.5, 3.1) -- (4.5, 3.1) ;
    \draw[gridColor] (-1.5, 4.1) -- (4.5, 4.1) ;
    \draw[gridColor, dashed] (-1, 4.1) -- (-1, 1.1) ;
    
    \draw[gridColor, dashed] (-3.5, -1.4) -- (-1, 1.1);
    \draw[gridColor, dashed] (-3.5, -0.4) -- (-1, 2.1) ;
    \draw[gridColor, dashed] (-3.5, 0.6) -- (-1, 3.1) ;
    \draw[gridColor] (-3.5, 1.6) -- (-1, 4.1) ;
    
    \draw[gridColor] (2, -1.4) -- (4.5, 1.1) ;
    \draw[gridColor] (2, -0.4) -- (4.5, 2.1) ;
    \draw[gridColor] (2, 0.6) -- (4.5, 3.1) ;
    \draw[gridColor] (2, 1.6) -- (4.5, 4.1) ;
    \draw[gridColor] (2, 1.6) -- (2, -1.4) ;
    \draw[gridColor] (4.5, 4.1) -- (4.5, 1.1) ;
    
    \tikzstyle wellColor=[plum1]
    \draw (0.5, 3)[wellColor, thick]  ellipse(0.5 and 0.2);
    \draw (0.0, 2)[wellColor, thick] arc[start angle=180,end angle=360,x radius=0.5,y radius=0.2]; 
    \draw (0.0, 1)[wellColor, thick] arc[start angle=180,end angle=360,x radius=0.5,y radius=0.2]; 
    \draw (0.0, 0)[wellColor, thick] arc[start angle=180,end angle=360,x radius=0.5,y radius=0.2]; 
    \draw[wellColor, thick] (0, 0) -- (0, 3) [] {};
    \draw[wellColor, thick] (1, 0) -- (1, 3) [] {};
    
    \node[] at (0.5, 0.4) [] {\footnotesize$p^\well_1$};
    \node[] at (0.5, 1.4) [] {\footnotesize$p^\well_2$};
    \node[] at (0.5, 2.4) [] {\footnotesize$p^\well_3$};
    
    \tikzstyle perfColor=[chocolate1]
    \draw[perfColor, thick] (-0.4, 0.35) -- (0.2, 0.35) [] {};
    \draw[perfColor, thick] (-0.4, 0.4) -- (0.2, 0.4) [] {};
    \draw[perfColor, thick] (-0.4, 1.35) -- (0.2, 1.35) [] {};
    \draw[perfColor, thick] (-0.4, 1.4) -- (0.2, 1.4) [] {};
    \draw[perfColor, thick] (-0.4, 2.35) -- (0.2, 2.35) [] {};
    \draw[perfColor, thick] (-0.4, 2.4) -- (0.2, 2.4) [] {};
    
    \node[] at (-0.25, 0.65) [] {\footnotesize $q^\well_1$};
    \node[] at (-0.25, 1.65) [] {\footnotesize $q^\well_2$};
    \node[] at (-0.25, 2.65) [] {\footnotesize $q^\well_3$};

    \node[] at (3.25, 0.35) [black!50] {\scriptsize $p^\res_{K(1)}$};
    \node[] at (3.25, 1.35) [black!50] {\scriptsize $p^\res_{K(2)}$};
    \node[] at (3.25, 2.35) [black!50] {\scriptsize $p^\res_{K(3)}$};
    
\end{tikzpicture}

%% file: pics/app_mesh_geometry.tex
\begin{tikzpicture}
[
scale=0.55,
cellCentroid/.style={circle,draw=black,text=black,fill=black,thick,minimum size=0.5mm, inner sep = 0em, outer sep = 0.0em},
faceCentroid/.style={rectangle,draw=black,text=black,fill=black,thick,minimum size=0.6mm, inner sep = 0em, outer sep = 0.0em},
cellLabel/.style={circle,draw=black,text=black,fill=red!50!black!20,minimum size=1mm, inner sep = 0.03em, outer sep = 0.1em},
faceLabel/.style={rectangle,draw=black,text=black,fill=blue!50!black!20,minimum size=2.5mm, inner sep = 0em, outer sep = 0.1em},
injectionWell/.style={circle,draw=red!50!black,fill=red!50!black!20,thick,minimum size=1mm, inner sep = 0.1em, outer sep = 0.15em},
productionWell/.style={circle,draw=blue!50!black,fill=blue!50!black!20,thick,minimum size=1mm, inner sep = 0.1em, outer sep = 0.15em}]

\coordinate [label=] (v1f) at (0, 0);
\coordinate [label=] (v2f) at (3, 0);
\coordinate [label=] (v3f) at (6, 0);
\coordinate [label=] (v4f) at (9, 0);
\coordinate [label=] (v5f) at (0, 2);
\coordinate [label=] (v6f) at (3, 2);
\coordinate [label=] (v7f) at (6, 2);
\coordinate [label=] (v8f) at (9, 2);
\coordinate [label=] (v9f) at (0, 4);
\coordinate [label=] (v10f) at (3, 4);
\coordinate [label=] (v11f) at (6, 4);
\coordinate [label=] (v12f) at (9, 4);

\coordinate [label=] (v1b) at (1, 1);
\coordinate [label=] (v2b) at (4, 1);
\coordinate [label=] (v3b) at (7, 1);
\coordinate [label=] (v4b) at (10, 1);
\coordinate [label=] (v5b) at (1, 3);
\coordinate [label=] (v6b) at (4, 3);
\coordinate [label=] (v7b) at (7, 3);
\coordinate [label=] (v8b) at (10, 3);
\coordinate [label=] (v9b) at (1, 5);
\coordinate [label=] (v10b) at (4, 5);
\coordinate [label=] (v11b) at (7, 5);
\coordinate [label=] (v12b) at (10, 5);

\tikzstyle gridColor=[black!30]

\draw[gridColor, thick] (v1f) -- (v4f);
\draw[gridColor, thick] (v5f) -- (v8f);
\draw[gridColor, thick] (v9f) -- (v12f);
\draw[gridColor, dashed] (v1b) -- (v4b);
\draw[gridColor, dashed] (v5b) -- (v8b);
\draw[gridColor, thick] (v9b) -- (v12b);

\draw[gridColor, thick] (v1f) -- (v9f);
\draw[gridColor, thick] (v2f) -- (v10f);
\draw[gridColor, thick] (v3f) -- (v11f);
\draw[gridColor, thick] (v4f) -- (v12f);
\draw[gridColor, dashed] (v1b) -- (v9b);
\draw[gridColor, dashed] (v2b) -- (v10b);
\draw[gridColor, dashed] (v3b) -- (v11b);
\draw[gridColor, thick] (v4b) -- (v12b);

\draw[gridColor, dashed] (v1f) -- (v1b) ;
\draw[gridColor, dashed] (v2f) -- (v2b) ;
\draw[gridColor, dashed] (v3f) -- (v3b) ;
\draw[gridColor, thick] (v4f) -- (v4b) ;
\draw[gridColor, dashed] (v5f) -- (v5b) ;
\draw[gridColor, dashed] (v6f) -- (v6b) ;
\draw[gridColor, dashed] (v7f) -- (v7b) ;
\draw[gridColor, thick] (v8f) -- (v8b) ;
\draw[gridColor, thick] (v9f) -- (v9b) ;
\draw[gridColor, thick] (v10f) -- (v10b) ;
\draw[gridColor, thick] (v11f) -- (v11b) ;
\draw[gridColor, thick] (v12f) -- (v12b) ;

\node[] at (barycentric cs:v1f=0.5,v6b=0.5) [cellLabel] {$\tau_1$};
\node[] at (barycentric cs:v6f=0.5,v9b=0.5) [cellLabel] {$\tau_2$};
\node[] at (barycentric cs:v2f=0.5,v7b=0.5) [cellLabel] {$\tau_3$};
\node[] at (barycentric cs:v6f=0.5,v11b=0.5) [cellLabel] {$\tau_4$};
\node[] at (barycentric cs:v4f=0.5,v7b=0.5) [cellLabel] {$\tau_5$};
\node[] at (barycentric cs:v7f=0.5,v12b=0.5) [cellLabel] {$\tau_6$};

\node[] at (barycentric cs:v2f=0.5,v6b=0.5) [faceCentroid] {};
\node[] at (barycentric cs:v5f=0.5,v6b=0.5) [faceCentroid] {};
\node[] at (barycentric cs:v3f=0.5,v7b=0.5) [faceCentroid] {};
\node[] at (barycentric cs:v6f=0.5,v7b=0.5) [faceCentroid] {};
\node[] at (barycentric cs:v7f=0.5,v8b=0.5) [faceCentroid] {};
\node[] at (barycentric cs:v6f=0.5,v10b=0.5) [faceCentroid] {};
\node[] at (barycentric cs:v7f=0.5,v11b=0.5) [faceCentroid] {};

\draw[thick,-stealth] (barycentric cs:v2f=0.5,v6b=0.5) -- node[above] {$\tensorOne{n}_{2}$} ++(0.1\linewidth,0);
\draw[thick,-stealth] (barycentric cs:v5f=0.5,v6b=0.5) -- node[right] {$\tensorOne{n}_{1}$} ++(0,0.1\linewidth);
\draw[thick,-stealth] (barycentric cs:v3f=0.5,v7b=0.5) -- node[above] {$\tensorOne{n}_{5}$} ++(0.1\linewidth,0);
\draw[thick,-stealth] (barycentric cs:v6f=0.5,v7b=0.5) -- node[right] {$\tensorOne{n}_{4}$} ++(0,0.1\linewidth);
\draw[thick,-stealth] (barycentric cs:v7f=0.5,v8b=0.5) -- node[right] {$\tensorOne{n}_{7}$} ++(0,0.1\linewidth);
\draw[thick,-stealth] (barycentric cs:v6f=0.5,v10b=0.5) -- node[above] {$\tensorOne{n}_{3}$} ++(0.1\linewidth,0);
\draw[thick,-stealth] (barycentric cs:v7f=0.5,v11b=0.5) -- node[above] {$\tensorOne{n}_{6}$} ++(0.1\linewidth,0);
\node[] at (barycentric cs:v5f=0.62,v6b=0.38) [faceLabel] {$\varepsilon_1$};
\node[] at (barycentric cs:v6f=0.62,v7b=0.38) [faceLabel] {$\varepsilon_4$};
\node[] at (barycentric cs:v7f=0.62,v8b=0.38) [faceLabel] {$\varepsilon_7$};
\node[] at (barycentric cs:v2f=0.68,v6b=0.38) [faceLabel] {$\varepsilon_2$};
\node[] at (barycentric cs:v6f=0.68,v10b=0.38) [faceLabel] {$\varepsilon_3$};
\node[] at (barycentric cs:v3f=0.68,v7b=0.38) [faceLabel] {$\varepsilon_5$};
\node[] at (barycentric cs:v7f=0.68,v11b=0.38) [faceLabel] {$\varepsilon_6$};

\end{tikzpicture}

%% file: pics/app_mesh_dofs.tex
\begin{tikzpicture}
[
scale=0.55,
cellCentroid/.style={circle,draw=black,text=black,fill=black,thick,minimum size=0.8mm, inner sep = 0em, outer sep = 0.0em},
faceCentroid/.style={rectangle,draw=black,text=black,fill=black,thick,minimum size=0.6mm, inner sep = 0em, outer sep = 0.0em},
cellLabel/.style={circle,draw=black,text=black,fill=red!50!black!20,thick,minimum size=4mm, inner sep = 0.1em, outer sep = 0.1em},
faceLabel/.style={rectangle,draw=black,text=black,fill=blue!50!black!20,thick,minimum size=4mm, inner sep = 0em, outer sep = 0.1em},
injectionWell/.style={circle,draw=red!50!black,fill=red!50!black!20,thick,minimum size=2mm, inner sep = 0.1em, outer sep = 0.15em},
productionWell/.style={circle,draw=blue!50!black,fill=blue!50!black!20,thick,minimum size=2mm, inner sep = 0.1em, outer sep = 0.15em}]

\coordinate [label=] (v1f) at (0, 0);
\coordinate [label=] (v2f) at (3, 0);
\coordinate [label=] (v3f) at (6, 0);
\coordinate [label=] (v4f) at (9, 0);
\coordinate [label=] (v5f) at (0, 2);
\coordinate [label=] (v6f) at (3, 2);
\coordinate [label=] (v7f) at (6, 2);
\coordinate [label=] (v8f) at (9, 2);
\coordinate [label=] (v9f) at (0, 4);
\coordinate [label=] (v10f) at (3, 4);
\coordinate [label=] (v11f) at (6, 4);
\coordinate [label=] (v12f) at (9, 4);

\coordinate [label=] (v1b) at (1, 1);
\coordinate [label=] (v2b) at (4, 1);
\coordinate [label=] (v3b) at (7, 1);
\coordinate [label=] (v4b) at (10, 1);
\coordinate [label=] (v5b) at (1, 3);
\coordinate [label=] (v6b) at (4, 3);
\coordinate [label=] (v7b) at (7, 3);
\coordinate [label=] (v8b) at (10, 3);
\coordinate [label=] (v9b) at (1, 5);
\coordinate [label=] (v10b) at (4, 5);
\coordinate [label=] (v11b) at (7, 5);
\coordinate [label=] (v12b) at (10, 5);

\tikzstyle gridColor=[black!15]

\draw[gridColor, thick] (v1f) -- (v4f);
\draw[gridColor, thick] (v5f) -- (v8f);
\draw[gridColor, thick] (v9f) -- (v12f);
\draw[gridColor, dashed] (v1b) -- (v4b);
\draw[gridColor, dashed] (v5b) -- (v8b);
\draw[gridColor, thick] (v9b) -- (v12b);

\draw[gridColor, thick] (v1f) -- (v9f);
\draw[gridColor, thick] (v2f) -- (v10f);
\draw[gridColor, thick] (v3f) -- (v11f);
\draw[gridColor, thick] (v4f) -- (v12f);
\draw[gridColor, dashed] (v1b) -- (v9b);
\draw[gridColor, dashed] (v2b) -- (v10b);
\draw[gridColor, dashed] (v3b) -- (v11b);
\draw[gridColor, thick] (v4b) -- (v12b);

\draw[gridColor, dashed] (v1f) -- (v1b) ;
\draw[gridColor, dashed] (v2f) -- (v2b) ;
\draw[gridColor, dashed] (v3f) -- (v3b) ;
\draw[gridColor, thick] (v4f) -- (v4b) ;
\draw[gridColor, dashed] (v5f) -- (v5b) ;
\draw[gridColor, dashed] (v6f) -- (v6b) ;
\draw[gridColor, dashed] (v7f) -- (v7b) ;
\draw[gridColor, thick] (v8f) -- (v8b) ;
\draw[gridColor, thick] (v9f) -- (v9b) ;
\draw[gridColor, thick] (v10f) -- (v10b) ;
\draw[gridColor, thick] (v11f) -- (v11b) ;
\draw[gridColor, thick] (v12f) -- (v12b) ;
%
%
%
\node[label={south west:$\{{\color{red!50!black}p_1}$,${\color{red!50!black}s_1}\}$}] at (barycentric cs:v1f=0.5,v6b=0.5) [cellCentroid] {};
\node[label={south:$\{{\color{red!50!black}p_3}$,${\color{red!50!black}s_3}\}$}] at (barycentric cs:v2f=0.5,v7b=0.5) [cellCentroid] {};
\node[label={south east:$\{{\color{red!50!black}p_5}$,${\color{red!50!black}s_5}\}$}] at (barycentric cs:v4f=0.5,v7b=0.5) [cellCentroid] {};
\node[label={west:$\{{\color{red!50!black}p_2}$,${\color{red!50!black}s_2}\}$}] at (barycentric cs:v6f=0.5,v9b=0.5) [cellCentroid] {};
\node[label={north:$\{{\color{red!50!black}p_4}$,${\color{red!50!black}s_4}\}$}] at (barycentric cs:v6f=0.5,v11b=0.5) [cellCentroid] {};
\node[label={east:$\{{\color{red!50!black}p_6}$,${\color{red!50!black}s_6}\}$}] at (barycentric cs:v7f=0.5,v12b=0.5) [cellCentroid] {};

\node[label={east,text=blue!50!black:$\sigma_1$}] at (barycentric cs:v5f=0.5,v6b=0.5)  [faceCentroid] {};
\node[label={east,text=blue!50!black:$\sigma_4$}] at (barycentric cs:v6f=0.5,v7b=0.5)  [faceCentroid] {};
\node[label={east,text=blue!50!black:$\sigma_7$}] at (barycentric cs:v7f=0.5,v8b=0.5)  [faceCentroid] {};
\node[label={south,text=blue!50!black:$\sigma_2$}]  at (barycentric cs:v2f=0.5,v6b=0.5)  [faceCentroid] {};
\node[label={south,text=blue!50!black:$\sigma_5$}]  at (barycentric cs:v3f=0.5,v7b=0.5)  [faceCentroid] {};
\node[label={south,text=blue!50!black:$\sigma_3$}]  at (barycentric cs:v6f=0.5,v10b=0.5) [faceCentroid] {};
\node[label={south,text=blue!50!black:$\sigma_6$}]  at (barycentric cs:v7f=0.5,v11b=0.5) [faceCentroid] {};

\end{tikzpicture}

%% file: pics/app_mesh_wells.tex
\begin{tikzpicture}
[
scale=0.55,
cellCentroid/.style={circle,draw=black,text=black,fill=black,minimum size=0.5mm, inner sep = 0em, outer sep = 0.0em},
perforationCentroid/.style={diamond,draw=black,text=black,fill=black,minimum size=0.5mm, inner sep = 0em, outer sep = 0.0em},
injectorLabel/.style={circle,draw=black,text=black,fill=green!50!black!20,minimum size=2mm, inner sep = 0.05em, outer sep = 0.1em},
producerLabel/.style={circle,draw=black,text=black,fill=yellow!50!black!20,minimum size=2mm, inner sep = 0.05em, outer sep = 0.1em},
perforationLabel/.style={diamond,draw=black,text=black,fill=blue!50!black!20,thick,minimum size=4mm, inner sep = 0em, outer sep = 0.1em}
]

\coordinate [label=] (v1f) at (0, 0);
\coordinate [label=] (v2f) at (3, 0);
\coordinate [label=] (v3f) at (6, 0);
\coordinate [label=] (v4f) at (9, 0);
\coordinate [label=] (v5f) at (0, 2);
\coordinate [label=] (v6f) at (3, 2);
\coordinate [label=] (v7f) at (6, 2);
\coordinate [label=] (v8f) at (9, 2);
\coordinate [label=] (v9f) at (0, 4);
\coordinate [label=] (v10f) at (3, 4);
\coordinate [label=] (v11f) at (6, 4);
\coordinate [label=] (v12f) at (9, 4);

\coordinate [label=] (v1b) at (1, 1);
\coordinate [label=] (v2b) at (4, 1);
\coordinate [label=] (v3b) at (7, 1);
\coordinate [label=] (v4b) at (10, 1);
\coordinate [label=] (v5b) at (1, 3);
\coordinate [label=] (v6b) at (4, 3);
\coordinate [label=] (v7b) at (7, 3);
\coordinate [label=] (v8b) at (10, 3);
\coordinate [label=] (v9b) at (1, 5);
\coordinate [label=] (v10b) at (4, 5);
\coordinate [label=] (v11b) at (7, 5);
\coordinate [label=] (v12b) at (10, 5);

\tikzstyle gridColor=[black!30]

\draw[gridColor, dashed] (v1b) -- (v4b);
\draw[gridColor, dashed] (v5b) -- (v8b);
\draw[gridColor, thick] (v9b) -- (v12b);

\draw[black!5!yellow, line width=3pt, line cap=round] (barycentric cs:v1f=0.5,v2b=0.5) -- ++(0,4.75);
\draw[black!30!green, line width=3pt, line cap=round] (barycentric cs:v3f=0.5,v4b=0.5) -- ++(0,4.75);

\draw[gridColor, thick] (v1f) -- (v4f);
\draw[gridColor, thick] (v5f) -- (v8f);
\draw[gridColor, thick] (v9f) -- (v12f);

\draw[gridColor, thick] (v1f) -- (v9f);
\draw[gridColor, thick] (v2f) -- (v10f);
\draw[gridColor, thick] (v3f) -- (v11f);
\draw[gridColor, thick] (v4f) -- (v12f);
\draw[gridColor, dashed] (v1b) -- (v9b);
\draw[gridColor, dashed] (v2b) -- (v10b);
\draw[gridColor, dashed] (v3b) -- (v11b);
\draw[gridColor, thick] (v4b) -- (v12b);

\draw[gridColor, dashed] (v1f) -- (v1b) ;
\draw[gridColor, dashed] (v2f) -- (v2b) ;
\draw[gridColor, dashed] (v3f) -- (v3b) ;
\draw[gridColor, thick] (v4f) -- (v4b) ;
\draw[gridColor, dashed] (v5f) -- (v5b) ;
\draw[gridColor, dashed] (v6f) -- (v6b) ;
\draw[gridColor, dashed] (v7f) -- (v7b) ;
\draw[gridColor, thick] (v8f) -- (v8b) ;
\draw[gridColor, thick] (v9f) -- (v9b) ;
\draw[gridColor, thick] (v10f) -- (v10b) ;
\draw[gridColor, thick] (v11f) -- (v11b) ;
\draw[gridColor, thick] (v12f) -- (v12b) ;

%
%
%

\node[] at (barycentric cs:v1f=0.5,v6b=0.5) [perforationCentroid] {};
\node[] at (barycentric cs:v3f=0.5,v8b=0.5) [perforationCentroid] {};
\node[] at (barycentric cs:v5f=0.5,v10b=0.5) [perforationCentroid] {};
\node[] at (barycentric cs:v7f=0.5,v12b=0.5) [perforationCentroid] {};


\draw[thick,-stealth] (barycentric cs:v1f=0.5,v6b=0.5) -- node[left] {} ++(-0.1\linewidth,0);
\draw[thick,-stealth] (barycentric cs:v1f=0.5,v6b=0.5) -- node[left] {} ++(0.1\linewidth,0);
\draw[thick,-stealth] (barycentric cs:v1f=0.5,v6b=0.5) -- node[left] {} ++(-0.09\linewidth,-0.09\linewidth);
\draw[thick,-stealth] (barycentric cs:v1f=0.5,v6b=0.5) -- node[left] {} ++(0.09\linewidth,0.09\linewidth);

\draw[thick,-stealth] (barycentric cs:v3f=0.5,v8b=0.5) -- node[left] {} ++(-0.1\linewidth,0);
\draw[thick,-stealth] (barycentric cs:v3f=0.5,v8b=0.5) -- node[left] {} ++(0.1\linewidth,0);
\draw[thick,-stealth] (barycentric cs:v3f=0.5,v8b=0.5) -- node[left] {} ++(-0.09\linewidth,-0.09\linewidth);
\draw[thick,-stealth] (barycentric cs:v3f=0.5,v8b=0.5) -- node[left] {} ++(0.09\linewidth,0.09\linewidth);

\draw[thick,-stealth] (barycentric cs:v5f=0.5,v10b=0.5) -- node[left] {} ++(-0.1\linewidth,0);
\draw[thick,-stealth] (barycentric cs:v5f=0.5,v10b=0.5) -- node[left] {} ++(0.1\linewidth,0);
\draw[thick,-stealth] (barycentric cs:v5f=0.5,v10b=0.5) -- node[left] {} ++(-0.09\linewidth,-0.09\linewidth);
\draw[thick,-stealth] (barycentric cs:v5f=0.5,v10b=0.5) -- node[left] {} ++(0.09\linewidth,0.09\linewidth);

\draw[thick,-stealth] (barycentric cs:v7f=0.5,v12b=0.5) -- node[left] {} ++(-0.1\linewidth,0);
\draw[thick,-stealth] (barycentric cs:v7f=0.5,v12b=0.5) -- node[left] {} ++(0.1\linewidth,0);
\draw[thick,-stealth] (barycentric cs:v7f=0.5,v12b=0.5) -- node[left] {} ++(-0.09\linewidth,-0.09\linewidth);
\draw[thick,-stealth] (barycentric cs:v7f=0.5,v12b=0.5) -- node[left] {} ++(0.09\linewidth,0.09\linewidth);

\draw[thick,-stealth] (1.5, 5.1) -- (2.6, 5.9) [] {};
\draw[thick,-stealth] (7.5, 5.8) -- (8.6, 5.1) [] {};

\node[label={north:\shortstack{BHP-controlled $(p^{\well, target}_1)$}}] at (2, 5.7) [] {};
\node[label={north:\shortstack{Rate-controlled $(q^{\well, target}_2)$}}] at (8, 5.7) [] {};
\node[] at (2, 5.5) [producerLabel] {$\well_1$};
\node[] at (8, 5.5) [injectorLabel] {$\well_2$};

\end{tikzpicture}

%% file: pics/app_well_dofs.tex
\begin{tikzpicture}
[
scale=0.55,
wellCentroid/.style={diamond,draw=black,text=black,fill=black,thick,minimum size=1mm, inner sep = 0em, outer sep = 0.0em},
perforationCentroid/.style={star,star points=6,draw=black,text=black,fill=black,thick,minimum size=1mm, inner sep = 0em, outer sep = 0.0em},
cellLabel/.style={circle,draw=black,text=black,fill=red!50!black!20,thick,minimum size=4mm, inner sep = 0.1em, outer sep = 0.1em},
perforationLabel/.style={diamond,draw=black,text=black,fill=blue!50!black!20,thick,minimum size=4mm, inner sep = 0em, outer sep = 0.1em},
injectionWell/.style={circle,draw=red!50!black,fill=red!50!black!20,thick,minimum size=2mm, inner sep = 0.1em, outer sep = 0.15em},
productionWell/.style={circle,draw=blue!50!black,fill=blue!50!black!20,thick,minimum size=2mm, inner sep = 0.1em, outer sep = 0.15em}]

\coordinate [label=] (v1f) at (0, 0);
\coordinate [label=] (v2f) at (3, 0);
\coordinate [label=] (v3f) at (6, 0);
\coordinate [label=] (v4f) at (9, 0);
\coordinate [label=] (v5f) at (0, 2);
\coordinate [label=] (v6f) at (3, 2);
\coordinate [label=] (v7f) at (6, 2);
\coordinate [label=] (v8f) at (9, 2);
\coordinate [label=] (v9f) at (0, 4);
\coordinate [label=] (v10f) at (3, 4);
\coordinate [label=] (v11f) at (6, 4);
\coordinate [label=] (v12f) at (9, 4);

\coordinate [label=] (v1b) at (1, 1);
\coordinate [label=] (v2b) at (4, 1);
\coordinate [label=] (v3b) at (7, 1);
\coordinate [label=] (v4b) at (10, 1);
\coordinate [label=] (v5b) at (1, 3);
\coordinate [label=] (v6b) at (4, 3);
\coordinate [label=] (v7b) at (7, 3);
\coordinate [label=] (v8b) at (10, 3);
\coordinate [label=] (v9b) at (1, 5);
\coordinate [label=] (v10b) at (4, 5);
\coordinate [label=] (v11b) at (7, 5);
\coordinate [label=] (v12b) at (10, 5);

\tikzstyle gridColor=[black!15]

\draw[gridColor, dashed] (v1b) -- (v4b);
\draw[gridColor, dashed] (v5b) -- (v8b);
\draw[gridColor, thick] (v9b) -- (v12b);

\draw[black!5!yellow, line width=3pt, line cap=round] (barycentric cs:v1f=0.5,v2b=0.5) -- ++(0,4.85);
\draw[black!30!green, line width=3pt, line cap=round] (barycentric cs:v3f=0.5,v4b=0.5) -- ++(0,4.85);

\draw[gridColor, thick] (v1f) -- (v4f);
\draw[gridColor, thick] (v5f) -- (v8f);
\draw[gridColor, thick] (v9f) -- (v12f);

\draw[gridColor, thick] (v1f) -- (v9f);
\draw[gridColor, thick] (v2f) -- (v10f);
\draw[gridColor, thick] (v3f) -- (v11f);
\draw[gridColor, thick] (v4f) -- (v12f);
\draw[gridColor, dashed] (v1b) -- (v9b);
\draw[gridColor, dashed] (v2b) -- (v10b);
\draw[gridColor, dashed] (v3b) -- (v11b);
\draw[gridColor, thick] (v4b) -- (v12b);

\draw[gridColor, dashed] (v1f) -- (v1b) ;
\draw[gridColor, dashed] (v2f) -- (v2b) ;
\draw[gridColor, dashed] (v3f) -- (v3b) ;
\draw[gridColor, thick] (v4f) -- (v4b) ;
\draw[gridColor, dashed] (v5f) -- (v5b) ;
\draw[gridColor, dashed] (v6f) -- (v6b) ;
\draw[gridColor, dashed] (v7f) -- (v7b) ;
\draw[gridColor, thick] (v8f) -- (v8b) ;
\draw[gridColor, thick] (v9f) -- (v9b) ;
\draw[gridColor, thick] (v10f) -- (v10b) ;
\draw[gridColor, thick] (v11f) -- (v11b) ;
\draw[gridColor, thick] (v12f) -- (v12b) ;

\node[label={west:${\color{red!50!black}\sigma_1^\well}$}] at (barycentric cs:v1f=0.5,v6b=0.5) [perforationCentroid] {};
\node[label={west:${\color{red!50!black}\sigma_5^\well}$}] at (barycentric cs:v3f=0.5,v8b=0.5) [perforationCentroid] {};
\node[label={west:${\color{red!50!black}\sigma_2^\well}$}] at (barycentric cs:v5f=0.5,v10b=0.5) [perforationCentroid] {};
\node[label={west:${\color{red!50!black}\sigma_6^\well}$}] at (barycentric cs:v7f=0.5,v12b=0.5) [perforationCentroid] {};

\node[label={west:$\color{green!30!black}p_2^\well$}] at (8, 5.35) [wellCentroid] {};

\end{tikzpicture}

%% file: pics/app_mat_M.tex
\begin{tikzpicture}
  \matrix[
    matrix of nodes,
    text height=2ex,
    text depth=0.0ex,
    text width=2ex,
    align=center,
    nodes={draw=black!10}, 
    nodes in empty cells,
  ] at (0,0) (M){
    &&&&&&&&&&\\
    &&&&&&&&&&\\
    &&&&&&&&&&\\
    &&&&&&&&&&\\
    &&&&&&&&&&\\
    &&&&&&&&&&\\
    &&&&&&&&&&\\
    &&&&&&&&&&\\
    &&&&&&&&&&\\
    &&&&&&&&&&\\
    &&&&&&&&&&\\
  };
  \path[fill=blue!50!black!20]
       (M-2-2.north west)
    -- (M-2-2.south west)
    -- (M-2-2.south east)
    -- (M-2-2.north east)
    -- cycle;  
  \path[fill=black!5!yellow]
       (M-9-9.north west)
    -- (M-9-9.south west)
    -- (M-9-9.south east)
    -- (M-9-9.north east)
    -- cycle; 
  \draw[thick,dashed,draw]
       (M-1-7.north east)
    -- (M-11-7.south east);
  \draw[thick,dashed,draw]
       (M-7-1.south west)
    -- (M-7-11.south east);  
  \path[thick,draw]
       (M-1-1.north west)
    -- (M-11-1.south west)
    -- (M-11-11.south east)
    -- (M-1-11.north east)
    -- cycle; 
  \foreach \i in {1,...,11}
    {
      \node[] at (M-\i-\i){$\bullet$};
    } 
  \draw (M-2-2.east) edge[out=350,in=160,-stealth,thick] (M-3-4.north west);  \node[anchor=west,fill=blue!50!black!20,draw,thick, outer sep=0.1em] at (M-3-4.north west) {$\frac{1}{\lambda(s_1) \overline{\Upsilon}_{1,\varepsilon}} + \frac{1}{\lambda(s_3) \overline{\Upsilon}_{3,\varepsilon}}$};
  \draw (M-9-9.west) edge[out=180 ,in=360,-stealth,thick] (M-10-5.east);  \node[anchor=east,fill=black!5!yellow,draw,thick, outer sep=0.1em] at (M-10-5.east) {$\frac{1}{\lambda(s_2) WI_2}$};

  \foreach \j in {1,...,7}
  {
    \node[anchor=south] at (M-1-\j.north){ $\sigma_\j$};
  }   
  \node[anchor=south] at (M-1-8.north){ $\sigma_1^\well$};
  \node[anchor=south] at (M-1-9.north){ $\sigma_2^\well$};
  \node[anchor=south] at (M-1-10.north){ $\sigma_5^\well$};
  \node[anchor=south] at (M-1-11.north){ $\sigma_6^\well$};
  
  \foreach \i in {1,...,7}
  {
    \node[anchor=east] at (M-\i-1.west){ $[\Vec{r}_{\sigma^\res}]_{\i}$};
  }
  \node[anchor=east] at (M-8-1.west){ $[\Vec{r}_{\sigma^\well}]_{1}$};
  \node[anchor=east] at (M-9-1.west){ $[\Vec{r}_{\sigma^\well}]_{2}$};
  \node[anchor=east] at (M-10-1.west){ $[\Vec{r}_{\sigma^\well}]_{3}$};
  \node[anchor=east] at (M-11-1.west){ $[\Vec{r}_{\sigma^\well}]_{4}$};
\end{tikzpicture}  

%% file: pics/app_vec_g.tex
\begin{tikzpicture}
  
  \matrix[
    matrix of nodes,
    text height=2ex,
    text depth=0.0ex,
    text width=2ex,
    align=center,
    nodes={draw=black!10}, 
    nodes in empty cells,
  ] at (0,0) (gVec){
    \\
    \\
    \\
    \\
    \\
    \\
    \\
    \\
    \\
    \\
    \\
  };
  \path[fill=black!5!yellow]
       (gVec-9-1.north west)
    -- (gVec-9-1.south west)
    -- (gVec-9-1.south east)
    -- (gVec-9-1.north east)
    -- cycle; 
  \draw[thick,dashed,draw]
       (gVec-7-1.south west)
    -- (gVec-7-1.south east);  
  \path[thick,draw]
       (gVec-1-1.north west)
    -- (gVec-11-1.south west)
    -- (gVec-11-1.south east)
    -- (gVec-1-1.north east)
    -- cycle; 

  \node[] at (gVec-9-1){$\bullet$};

  \node[anchor=east,fill=black!5!yellow,draw,thick, outer sep=0.1em] (bhp) at (11em, -4em) {$-(p_{1}^{\well, target})$};
  \draw (gVec-9-1.east) edge[out=0,in=200,-stealth,thick] (bhp.west);
  

  \foreach \i in {1,...,7}
  {
    \node[anchor=east] at (gVec-\i-1.west){ $[\Vec{r}_{\sigma^\res}]_{\i}$};
  }
  \node[anchor=east] at (gVec-8-1.west){ $[\Vec{r}_{\sigma^\well}]_{1}$};
  \node[anchor=east] at (gVec-9-1.west){ $[\Vec{r}_{\sigma^\well}]_{2}$};
  \node[anchor=east] at (gVec-10-1.west){ $[\Vec{r}_{\sigma^\well}]_{3}$};
  \node[anchor=east] at (gVec-11-1.west){ $[\Vec{r}_{\sigma^\well}]_{4}$};
\end{tikzpicture} 

%% file: pics/app_mat_D.tex
\begin{tikzpicture}
  
  \matrix[
    matrix of nodes,
    text height=2ex,
    text depth=0.0ex,
    text width=2ex,
    align=center,
    nodes={draw=black!10}, 
    nodes in empty cells,
  ] at (0,0) (D){
    &&&&&&&&&&\\
    &&&&&&&&&&\\
    &&&&&&&&&&\\
    &&&&&&&&&&\\
    &&&&&&&&&&\\
    &&&&&&&&&&\\
    &&&&&&&&&&\\
  };
  \draw[thick,dashed,draw]
       (D-1-7.north east)
    -- (D-7-7.south east);
  \draw[thick,dashed,draw]
       (D-6-1.south west)
    -- (D-6-11.south east);
  \path[fill=black!30!green]
       (D-7-10.north west)
    -- (D-7-11.north east)
    -- (D-7-11.south east)
    -- (D-7-10.south west)
    -- cycle; 
  \path[fill=black!5!yellow]
       (D-2-9.north west)
    -- (D-2-9.south west)
    -- (D-2-9.south east)
    -- (D-2-9.north east)
    -- cycle; 
  \path[thick,draw]
       (D-1-1.north west)
    -- (D-1-11.north east)
    -- (D-7-11.south east)
    -- (D-7-1.south west)
    -- cycle; 

  \node[] at (D-1-1){1};
  \node[] at (D-2-1){-1};
  \node[] at (D-1-2){1};
  \node[] at (D-3-2){-1};
  \node[] at (D-2-3){1};
  \node[] at (D-4-3){-1};
  \node[] at (D-3-4){1};
  \node[] at (D-4-4){-1};
  \node[] at (D-3-5){1};
  \node[] at (D-5-5){-1};
  \node[] at (D-4-6){1};
  \node[] at (D-6-6){-1};
  \node[] at (D-5-7){1};
  \node[] at (D-6-7){-1};
  \node[] at (D-1-8){1};
  \node[] at (D-2-9){1};
  \node[] at (D-5-10){1};
  \node[] at (D-7-10){-1};
  \node[] at (D-6-11){1};
  \node[] at (D-7-11){-1};

  \foreach \i in {1,...,6}
  {
    \node[anchor=east] at (D-\i-1.west){ $[\Vec{r}_{p^\res}]_{\i}$};
  }
  \node[anchor=east] at (D-7-1.west){ $[\Vec{r}_{p^\well}]_1$};
  
  \foreach \i in {1,...,7}
  {
    \node[anchor=south] at (D-1-\i.north){ $\sigma^\res_\i$};
  }
  \node[anchor=south] at (D-1-8.north){ $\sigma^\well_1$};
  \node[anchor=south] at (D-1-9.north){ $\sigma^\well_2$};
  \node[anchor=south] at (D-1-10.north){ $\sigma^\well_5$};
  \node[anchor=south] at (D-1-11.north){ $\sigma^\well_6$};

\end{tikzpicture}  

%% file: pics/app_vec_f.tex
\begin{tikzpicture}
  
  \matrix[
    matrix of nodes,
    text height=2ex,
    text depth=0.0ex,
    text width=2ex,
    align=center,
    nodes={draw=black!10}, 
    nodes in empty cells,
  ] at (0,0) (gVec){
    \\
    \\
    \\
    \\
    \\
    \\
    \\
  };
  \path[fill=black!30!green]
       (gVec-7-1.north west)
    -- (gVec-7-1.south west)
    -- (gVec-7-1.south east)
    -- (gVec-7-1.north east)
    -- cycle; 
  \draw[thick,dashed,draw]
       (gVec-6-1.south west)
    -- (gVec-6-1.south east);  
  \path[thick,draw]
       (gVec-1-1.north west)
    -- (gVec-7-1.south west)
    -- (gVec-7-1.south east)
    -- (gVec-1-1.north east)
    -- cycle; 

  \node[] at (gVec-7-1){$\bullet$};

  \node[anchor=east,fill=black!30!green,draw,thick, outer sep=0.1em] (bhp) at (10.5em, -3em) {$q^{\well, target}_2$};
  \draw (gVec-7-1.east) edge[out=0,in=200,-stealth,thick] (bhp.west);
  
  \foreach \i in {1,...,6}
  {
    \node[anchor=east] at (gVec-\i-1.west){ $[\Vec{r}_{p^\res}]_{\i}$};
  }
  \node[anchor=east] at (gVec-7-1.west){ $[\Vec{r}_{p^\well}]_{1}$};
\end{tikzpicture} 

%% file: pics/well-aggregation.tex
\begin{tikzpicture}
[
scale=0.75
]
\tikzstyle gridColor=[green!50]
\tikzstyle edgeColor=[black!50]
\tikzstyle vertColor=[black!80]

\foreach \x in {1,2,...,11}
{
	\draw[gridColor, dashed, line width=0.5mm] (\x, 2) -- (\x, 11);
}
\foreach \x in {2,3,...,11}
{
	\draw[gridColor, dashed, line width=0.5mm] (1, \x) -- (11, \x);
}    

\foreach \x in {1.5,2.5,...,10.5}
{
	\draw[edgeColor] (\x, 2.5) -- (\x, 10.5);
}
\foreach \x in {2.5,3.5,...,10.5}
{
	\draw[edgeColor] (1.5, \x) -- (10.5, \x);
}

\tikzstyle weightColor=[red!70!black!50]
\draw[weightColor, line width=0.8mm] (5.5, 7.5) -- (9.5, 7.5);
\draw[weightColor, line width=0.8mm] (7.5, 5.5) -- (7.5, 9.5);
\draw[weightColor, line width=0.8mm] (6.5, 6.5) -- (6.5, 8.5);
\draw[weightColor, line width=0.8mm] (8.5, 6.5) -- (8.5, 8.5);
\draw[weightColor, line width=0.8mm] (6.5, 6.5) -- (8.5, 6.5);
\draw[weightColor, line width=0.8mm] (6.5, 8.5) -- (8.5, 8.5);

\foreach \x in {1.5,2.5,...,10.5}
{
	\foreach \y in {2.5,3.5,...,10.5}
	{
		\draw (\x, \y)[vertColor, thick, fill] circle (0.15);
	}
 }
   
\tikzstyle wellColor=[blue!70!black!50]
\draw (7.5, 7.5)[wellColor, thick, fill] circle (0.25);

\draw[gridColor, dashed, line width=0.5mm] (12.5, 8.5) -- (12.5, 9.5) -- (13.5, 9.5) -- (13.5, 8.5) -- (12.5, 8.5);
\node[anchor=west] at (14, 9) [] {: Mesh};
\draw[edgeColor] (12.5, 7) -- (12.5, 8) -- (13.5, 8) -- (13.5, 7) -- (12.5, 7);
\draw (12.5, 7)[vertColor, thick, fill] circle (0.15);
\draw (12.5, 8)[vertColor, thick, fill] circle (0.15);
\draw (13.5, 7)[vertColor, thick, fill] circle (0.15);
\draw (13.5, 8)[vertColor, thick, fill] circle (0.15);
\node[anchor=west] at (14, 7.5) [] {: cell-connectivity graph};
\draw (13, 6)[wellColor, thick, fill] circle (0.25);
\node[anchor=west] at (14, 6) [] {: well cell};
\draw[weightColor, line width=0.8mm] (13, 4) -- (13, 5);
\node[anchor=west] at (14, 4.5) [] {: edges with heavy weights};

\end{tikzpicture}

%% file: manuscript.bbl
\begin{thebibliography}{10}
\expandafter\ifx\csname url\endcsname\relax
  \def\url#1{\texttt{#1}}\fi
\expandafter\ifx\csname urlprefix\endcsname\relax\def\urlprefix{URL }\fi
\expandafter\ifx\csname href\endcsname\relax
  \def\href#1#2{#2} \def\path#1{#1}\fi

\bibitem{brandt77}
A.~Brandt, Multi-level adaptive solutions to boundary-value problems, Math.
  Comp. 31~(138) (1977) 333--390.
\newblock \href {https://doi.org/10.2307/2006422} {\path{doi:10.2307/2006422}}.

\bibitem{jenny2009}
P.~Jenny, H.~A. Tchelepi, S.~H. Lee, Unconditionally convergent nonlinear
  solver for hyperbolic conservation laws with {S}-shaped flux functions, J.
  Comput. Phys. 228~(20) (2009) 7497--7512.
\newblock \href {https://doi.org/10.1016/j.jcp.2009.06.032}
  {\path{doi:10.1016/j.jcp.2009.06.032}}.

\bibitem{wang2013}
X.~Wang, H.~A. Tchelepi, Trust-region based solver for nonlinear transport in
  heterogeneous porous media, J. Comput. Phys. 253 (2013) 114--137.
\newblock \href {https://doi.org/10.1016/j.jcp.2013.06.041}
  {\path{doi:10.1016/j.jcp.2013.06.041}}.

\bibitem{li2015}
B.~Li, H.~A. Tchelepi, Nonlinear analysis of multiphase transport in porous
  media in the presence of viscous, buoyancy, and capillary forces, J. Comput.
  Phys. 297 (2015) 104--131.
\newblock \href {https://doi.org/10.1016/j.jcp.2015.04.057}
  {\path{doi:10.1016/j.jcp.2015.04.057}}.

\bibitem{moyner2017}
O.~M{\o}yner, Nonlinear solver for three-phase transport problems based on
  approximate trust regions, Comput. Geosci. 21 (2017) 999--1021.
\newblock \href {https://doi.org/10.1007/s10596-017-9660-1}
  {\path{doi:10.1007/s10596-017-9660-1}}.

\bibitem{lee2015hybrid}
S.~H. Lee, Y.~Efendiev, H.~A. Tchelepi, Hybrid upwind discretization of
  nonlinear two-phase flow with gravity, Adv. Water Resour. 82 (2015) 27--38.
\newblock \href {https://doi.org/10.1016/j.advwatres.2015.04.007}
  {\path{doi:10.1016/j.advwatres.2015.04.007}}.

\bibitem{hamon2016implicit}
F.~P. Hamon, B.~Mallison, H.~A. Tchelepi, Implicit hybrid upwind scheme for
  coupled multiphase flow and transport with buoyancy, Comput. Methods Appl.
  Mech. Engrg. 311 (2016) 599--624.
\newblock \href {https://doi.org/10.1016/j.cma.2016.08.009}
  {\path{doi:10.1016/j.cma.2016.08.009}}.

\bibitem{moncorge2020consistent}
A.~Moncorg{\'e}, O.~M{\o}yner, H.~A. Tchelepi, P.~Jenny, Consistent upwinding
  for sequential fully implicit multiscale compositional simulation, Comput.
  Geosci. 24~(2) (2020) 533--550.
\newblock \href {https://doi.org/10.1007/s10596-019-09835-6}
  {\path{doi:10.1007/s10596-019-09835-6}}.

\bibitem{bosma2022smooth}
S.~Bosma, F.~P. Hamon, B.~Mallison, H.~A. Tchelepi, Smooth implicit hybrid
  upwinding for compositional multiphase flow in porous media, Comput. Methods
  Appl. Mech. Engrg. 388 (2022) 114288.
\newblock \href {https://doi.org/10.1016/j.cma.2021.114288}
  {\path{doi:10.1016/j.cma.2021.114288}}.

\bibitem{kwok2007}
F.~Kwok, H.~A. Tchelepi, Potential-based reduced {N}ewton algorithm for
  nonlinear multiphase flow in porous media, J. Comput. Phys. 227~(1) (2007)
  706--727.
\newblock \href {https://doi.org/10.1016/j.jcp.2007.08.012}
  {\path{doi:10.1016/j.jcp.2007.08.012}}.

\bibitem{natvig2008}
J.~Natvig, K.-A. Lie, Fast computation of multiphase flow in porous media by
  implicit discontinuous {G}alerkin schemes with optimal ordering of elements,
  J. Comput. Phys. 227~(24) (2008) 10108--10124.
\newblock \href {https://doi.org/10.1016/j.jcp.2008.08.024}
  {\path{doi:10.1016/j.jcp.2008.08.024}}.

\bibitem{hamon2016ordering}
F.~P. Hamon, H.~A. Tchelepi, Ordering-based nonlinear solver for fully implicit
  simulation of three-phase flow, Comput. Geosci. 20~(5) (2016) 909--927.
\newblock \href {https://doi.org/10.1007/s10596-016-9569-0}
  {\path{doi:10.1007/s10596-016-9569-0}}.

\bibitem{kelmetstal2020reordering}
{\O}.~S. Klemetsdal, A.~F. Rasmussen, O.~M{\o}yner, K.-A. Lie, Efficient
  reordered nonlinear {Gauss--Seidel} solvers with higher order for black-oil
  models, Comput. Geosci. 24~(2) (2020) 593--607.
\newblock \href {https://doi.org/10.1007/s10596-019-09844-5}
  {\path{doi:10.1007/s10596-019-09844-5}}.

\bibitem{cai2002}
X.-C. Cai, D.~E. Keyes, Nonlinearly preconditioned inexact {Newton} algorithms,
  SIAM J. Sci. Comput. 24~(1) (2002) 183--200.
\newblock \href {https://doi.org/10.1137/S106482750037620X}
  {\path{doi:10.1137/S106482750037620X}}.

\bibitem{liu2015}
L.~Liu, D.~E. Keyes, Field-split preconditioned inexact {Newton} algorithms,
  SIAM J. Sci. Comput. 37~(3) (2015) A1388--A1409.
\newblock \href {https://doi.org/10.1137/140970379}
  {\path{doi:10.1137/140970379}}.

\bibitem{dolean2016}
V.~Dolean, M.~J. Gander, W.~Kheriji, F.~Kwok, R.~Masson, Nonlinear
  preconditioning: How to use a nonlinear {Schwarz} method to precondition
  {Newton's} method, SIAM J. Sci. Comput. 38~(6) (2016) A3357--A3380.
\newblock \href {https://doi.org/10.1137/15M102887X}
  {\path{doi:10.1137/15M102887X}}.

\bibitem{skogestad2013}
J.~O. Skogestad, E.~Keilegavlen, J.~M. Nordbotten, Domain decomposition
  strategies for nonlinear flow problems in porous media, J. Comput. Phys. 234
  (2013) 439--451.
\newblock \href {https://doi.org/10.1016/j.jcp.2012.10.001}
  {\path{doi:10.1016/j.jcp.2012.10.001}}.

\bibitem{skogestad2016}
J.~O. Skogestad, E.~Keilegavlen, J.~M. Nordbotten, Two-scale preconditioning
  for two-phase nonlinear flows in porous media, Transp. Porous Media 114~(2)
  (2016) 485--503.
\newblock \href {https://doi.org/10.1007/s10596-019-09844-5}
  {\path{doi:10.1007/s10596-019-09844-5}}.

\bibitem{kelmetstal2020schwarz}
{\O}.~S. Klemetsdal, A.~Moncorg{\'e}, O.~M{\o}yner, K.-A. Lie, Additive
  {Schwarz} preconditioned exact newton method as a nonlinear preconditioner
  for multiphase porous media flow, in: Proceedings of ECMOR XVII, Vol. 2020,
  European Association of Geoscientists and Engineers, 2020, pp. 1--20.
\newblock \href {https://doi.org/10.3997/2214-4609.202035050}
  {\path{doi:10.3997/2214-4609.202035050}}.

\bibitem{n2023comparison}
M.~N'diaye, F.~P. Hamon, H.~A. Tchelepi, Comparison of nonlinear field-split
  preconditioners for two-phase flow in heterogeneous porous media, Comput.
  Geosci. (2023) 1--17\href {https://doi.org/10.1007/s10596-023-10200-x}
  {\path{doi:10.1007/s10596-023-10200-x}}.

\bibitem{christensen16}
M.~l.~C. Christensen, K.~L. Eskildsen, A.~P. Engsig-Karup, M.~A. Wakefield,
  Nonlinear multigrid for reservoir simulation, SPE J. 21~(3) (2016) 888--898.
\newblock \href {https://doi.org/10.2118/178428-PA}
  {\path{doi:10.2118/178428-PA}}.

\bibitem{christensen18}
M.~l.~C. Christensen, P.~S. Vassilevski, U.~Villa, Nonlinear multigrid solvers
  exploiting {AMG}e coarse spaces with approximation properties, J. Comput.
  Appl. Math. 340 (2018) 691--708.
\newblock \href {https://doi.org/10.1016/j.cam.2017.10.029}
  {\path{doi:10.1016/j.cam.2017.10.029}}.

\bibitem{toft18}
R.~Toft, K.-A. Lie, O.~M{\o}yner,
  \href{{https://ojs.bibsys.no/index.php/NIK/article/view/503}}{ Full
  approximation scheme for reservoir simulation}, in: Proceedings - Norsk
  Informatikkonferanse, Oslo, Norway, 2018.

\bibitem{fas-spectral-diffusion}
C.~S. Lee, F.~P. Hamon, N.~Castelletto, P.~S. Vassilevski, J.~A. White,
  Nonlinear multigrid based on local spectral coarsening for heterogeneous
  diffusion problems, Comput. Methods Appl. Mech. Engrg. 372 (2020) 113432.
\newblock \href {https://doi.org/10.1016/j.cma.2020.113432}
  {\path{doi:10.1016/j.cma.2020.113432}}.

\bibitem{fas-two-phase}
C.~S. Lee, F.~P. Hamon, N.~Castelletto, P.~S. Vassilevski, J.~A. White, An
  aggregation-based nonlinear multigrid solver for two-phase flow and transport
  in porous media, Comput. Math. with Appl. 113 (2022) 282--299.
\newblock \href {https://doi.org/https://doi.org/10.1016/j.camwa.2022.03.026}
  {\path{doi:https://doi.org/10.1016/j.camwa.2022.03.026}}.

\bibitem{Pea78}
D.~W. Peaceman, Interpretation of well-block pressures in numerical reservoir
  simulation, SPE J. 18~(3) (1978) 183--194.
\newblock \href {https://doi.org/10.2118/6893-PA} {\path{doi:10.2118/6893-PA}}.

\bibitem{wolfsteiner06}
C.~Wolfsteiner, S.~H. Lee, H.~A. Tchelepi, Well modeling in the multiscale
  finite volume method for subsurface flow simulation, Multiscale Model. Simul.
  5~(3) (2006) 900--917.
\newblock \href {https://doi.org/https://doi.org/10.1137/050640771}
  {\path{doi:https://doi.org/10.1137/050640771}}.

\bibitem{arbogast2002two}
T.~Arbogast, S.~L. Bryant, A two-scale numerical subgrid technique for
  waterflood simulations, SPE J. 7~(04) (2002) 446--457.
\newblock \href {https://doi.org/10.2118/81909-PA}
  {\path{doi:10.2118/81909-PA}}.

\bibitem{chen2003numerical}
Z.~Chen, X.~Yue, Numerical homogenization of well singularities in the flow
  transport through heterogeneous porous media, Multiscale Model. Simul. 1~(2)
  (2003) 260--303.
\newblock \href {https://doi.org/10.1137/S1540345902413322}
  {\path{doi:10.1137/S1540345902413322}}.

\bibitem{aarnes2004use}
J.~E. Aarnes, On the use of a mixed multiscale finite element method for
  greater flexibility and increased speed or improved accuracy in reservoir
  simulation, Multiscale Model. Simul. 2~(3) (2004) 421--439.
\newblock \href {https://doi.org/10.1137/030600655}
  {\path{doi:10.1137/030600655}}.

\bibitem{aarnes06}
J.~E. Aarnes, S.~Krogstad, K.-A. Lie, A hierarchical multiscale method for
  two-phase flow based upon mixed finite elements and nonuniform coarse grids,
  Multiscale Model. Simul. 5~(2) (2006) 337--363.
\newblock \href {https://doi.org/10.1137/050634566}
  {\path{doi:10.1137/050634566}}.

\bibitem{ligaarden-thesis}
S.~Ligaarden, Well models for mimetic finite difference methods and improved
  representation of wells in multiscale methods, Master's thesis, University of
  Oslo (2008).

\bibitem{skaflestad2008multiscale}
B.~Skaflestad, S.~Krogstad, Multiscale/mimetic pressure solvers with near-well
  grid adaptation, in: ECMOR XI-11th European Conference on the Mathematics of
  Oil Recovery, European Association of Geoscientists \& Engineers, 2008, pp.
  cp--62.
\newblock \href {https://doi.org/10.3997/2214-4609.20146387}
  {\path{doi:10.3997/2214-4609.20146387}}.

\bibitem{aziz79}
K.~Aziz, A.~Settari, Petroleum Reservoir Simulation, Elsevier Applied Science
  Publishers, London, UK, 1979.

\bibitem{EymGalHer00}
R.~Eymard, T.~Gallou{\"{e}}t, R.~Herbin, Finite volume methods, in: P.~G.
  Ciarlet, J.-L. Lions (Eds.), Solution of Equation in {$R^n$} (Part 3),
  Techniques of Scientific Computing (Part 3), Vol.~7 of Handbook of Numerical
  Analysis, Elsevier, 2000, pp. 713--1018.
\newblock \href {https://doi.org/10.1016/S1570-8659(00)07005-8}
  {\path{doi:10.1016/S1570-8659(00)07005-8}}.

\bibitem{lie19}
K.-A. Lie, An Introduction to Reservoir Simulation Using MATLAB/GNU Octave:
  User Guide for the MATLAB Reservoir Simulation Toolbox (MRST), Cambridge
  University Press, 2019.
\newblock \href {https://doi.org/10.1017/9781108591416}
  {\path{doi:10.1017/9781108591416}}.

\bibitem{CheHuaMa06}
Z.~Chen, G.~Huan, Y.~Ma, Computational Methods for Multiphase Flows in Porous
  Media, Society for Industrial and Applied Mathematics, 2006.
\newblock \href {https://doi.org/10.1137/1.9780898718942}
  {\path{doi:10.1137/1.9780898718942}}.

\bibitem{henson2003multigrid}
V.~E. Henson, Multigrid methods for nonlinear problems: {A}n overview, in:
  C.~A. Bouman, R.~L. Stevenson (Eds.), Computational Imaging, Vol. 5016 of
  Proceedings of SPIE, 2003, pp. 36--48.
\newblock \href {https://doi.org/10.1117/12.499473}
  {\path{doi:10.1117/12.499473}}.

\bibitem{karypis1998fast}
G.~Karypis, V.~Kumar, A fast and highly quality multilevel scheme for
  partitioning irregular graphs, SIAM J. Sci. Comput. 20~(1) (1998) 359--392.
\newblock \href {https://doi.org/10.1137/S1064827595287997}
  {\path{doi:10.1137/S1064827595287997}}.

\bibitem{ml-spectral-coarsening}
A.~T. Barker, S.~V. Gelever, C.~S. Lee, S.~Osborn, P.~S. Vassilevski,
  Multilevel spectral coarsening for graph {L}aplacian problems with
  application to reservoir simulation, SIAM J. Sci. Comput. 43~(4) (2021)
  A2737--A2765.
\newblock \href {https://doi.org/10.1137/19M1296343}
  {\path{doi:10.1137/19M1296343}}.

\bibitem{lee17}
C.~S. Lee, P.~S. Vassilevski, Parallel solver for {${H}$}(div) problems using
  hybridization and {AMG}, in: C.-O. Lee, X.-C. Cai, D.~E. Keyes, H.~H. Kim,
  A.~Klawonn, E.-J. Park, O.~B. Widlund (Eds.), Domain Decomposition Methods in
  Science and Engineering XXIIII, Vol. 116 of Lecture Notes in Computational
  Science and Engineering, 2017, pp. 69--80.
\newblock \href {https://doi.org/10.1007/978-3-319-52389-7_6}
  {\path{doi:10.1007/978-3-319-52389-7_6}}.

\bibitem{dobrev18}
V.~A. Dobrev, {\relax Tz}.~V. Kolev, C.~S. Lee, V.~Z. Tomov, P.~S. Vassilevski,
  Algebraic hybridization and static condensation with application to scalable
  {$H$(div)} preconditioning, SIAM J. Sci. Comput. 41~(3) (2019) B425--B447.
\newblock \href {https://doi.org/10.1137/17M1132562}
  {\path{doi:10.1137/17M1132562}}.

\bibitem{cao2002development}
H.~Cao, Development of techniques for general purpose simulators, Stanford
  University, 2002.

\bibitem{younis2011modern}
R.~Younis, Modern advances in software and solution algorithms for reservoir
  simulation, Stanford University, 2011.

\bibitem{mfem}
{MFEM}: {M}odular finite element methods library, \url{http://mfem.org}.

\bibitem{smoothg}
{smoothG}: {M}ixed graph {L}aplacian upscaling and solvers,
  \url{https://github.com/LLNL/smoothG}.

\bibitem{EggModel}
J.~D. Jansen, R.-M. Fonseca, S.~Kahrobaei, M.~M. Siraj, G.~M. Van~Essen,
  P.~M.~J. Van~den Hof, The egg model - a geological ensemble for reservoir
  simulation, Geosci. Data J. 1~(2) (2014) 192--195.
\newblock \href {https://doi.org/10.1002/gdj3.21} {\path{doi:10.1002/gdj3.21}}.

\bibitem{GEOS}
{GEOS}, \url{www.geosx.org}, [Online; accessed 2023-02-03] (2023).

\bibitem{manzocchi2008sensitivity}
T.~Manzocchi, J.~N. Carter, A.~Skorstad, B.~Fjellvoll, K.~D. Stephen, J.~A.
  Howell, J.~D. Matthews, J.~J. Walsh, M.~Nepveu, C.~Bos, Sensitivity of the
  impact of geological uncertainty on production from faulted and unfaulted
  shallow-marine oil reservoirs: objectives and methods, Petroleum Geoscience
  14~(1) (2008) 3--15.
\newblock \href {https://doi.org/10.1144/1354-079307-790}
  {\path{doi:10.1144/1354-079307-790}}.

\end{thebibliography}
